\documentclass{article}

\usepackage{pb-diagram}

\usepackage{amsmath}
\usepackage{amsthm}
\usepackage{amssymb}

\usepackage{bm}

\usepackage{latexsym}

  \newtheorem{theorem}{Theorem}[section]
  \newtheorem{corollary}[theorem]{Corollary}
  \newtheorem{lemma}[theorem]{Lemma}
  \newtheorem{proposition}[theorem]{Proposition}
  \newtheorem*{introtheorem}{Theorem}

  \newtheorem{definition}[theorem]{Definition}
  \newtheorem{convention}[theorem]{Convention}
  \newtheorem{ex}[theorem]{Example}
  \newtheorem{remark}[theorem]{Remark}

  \newtheorem*{acknowledgement}{Acknowledgements}

  \newenvironment{example}{\begin{ex}}{\qed\end{ex}}

\newcommand{\Sal}{\operatorname{Sal}} 
\newcommand{\R}{\mathbb{R}}
\newcommand{\C}{\mathbb{C}}
\newcommand{\Z}{\mathbb{Z}}
\newcommand{\Tot}{\operatorname{Tot}}
\newcommand{\Cotor}{\operatorname{Cotor}}
\newcommand{\Map}{\operatorname{Map}}
\newcommand{\Ima}{\operatorname{Im}}
\DeclareMathOperator*{\holim}{\mathrm{holim}}
\newcommand{\rarrow}[1]{\buildrel #1 \over \longrightarrow}
\newcommand{\larrow}[1]{\buildrel #1 \over \longleftarrow}
\newcommand{\shot}{\mathop\simeq\limits_S}
\newcommand{\sgn}{\operatorname{\mathrm{sgn}}}
\newcommand{\sign}{\operatorname{\mathrm{sgn}}}
\newcommand{\ad}{\operatorname{\mathrm{ad}}} 
\newcommand{\op}{{\operatorname{\mathit{op}}}} 
\newcommand{\sk}{\operatorname{sk}}
\newcommand{\Conv}{\operatorname{\mathrm{Conv}}} 
\newcommand{\Int}{{\operatorname{\mathrm{Int}}}}
\newcommand{\smallfrac}[2]{{\textstyle \frac{#1}{#2}}}
\newcommand{\codim}{\operatorname{\mathrm{codim}}}
\newcommand{\lindep}{\operatorname{lin-dep}} 
\newcommand{\alg}{\operatorname{\mathrm{alg}}}
\newcommand{\geo}{\operatorname{\mathrm{geo}}}
\newcommand{\lset}[2]{%
\left.\left\{#1 \ \right| \ #2\right\}
}

\newcommand{\bibdir}{bib}

\usepackage{hyperref}

\title{\bfseries The Salvetti Complex and the Little Cubes}
\author{Dai Tamaki
\thanks{Partially supported by Grants-in-Aid for Scientific Research,
Ministry of Education, Culture, Sports, Science and Technology, Japan:
17540070} \\
Department of Mathematical Sciences, \\
Shinshu University, \\
Matsumoto, 390-8621, Japan \\
\texttt{rivulus@math.shinshu-u.ac.jp}
}
\date{}

\begin{document}

\maketitle

\begin{abstract}
 For a real central arrangement $\mathcal{A}$, Salvetti introduced a
 construction of a finite complex $\Sal(\mathcal{A})$ which is homotopy
 equivalent to the complement of the complexified arrangement in
 \cite{Salvetti87}.  For the braid arrangement $\mathcal{A}_{k-1}$, the
 Salvetti complex $\Sal(\mathcal{A}_{k-1})$ serves as a good combinatorial
 model for the homotopy type of the configuration space $F(\C,k)$ of $k$
 points in $\C$, which is homotopy equivalent to the space
 $\mathcal{C}_2(k)$ of $k$ little $2$-cubes.
 Motivated by the importance of little cubes in homotopy theory,
 especially in the study of iterated loop spaces, we study how the
 combinatorial structure of the Salvetti complexes of the braid
 arrangements are related to homotopy theoretic properties of iterated
 loop spaces.

 As a consequence, we prove the skeletal filtrations on the Salvetti
 complexes of the braid arrangements give rise to the cobar-type
 Eilenberg-Moore spectral sequence converging to the homology of
 $\Omega^2\Sigma^2 X$. We also construct a new spectral sequence that
 computes the homology of $\Omega^{\ell}\Sigma^{\ell} X$ for $\ell>2$ by
 using a higher order analogue of the Salvetti complex. The
 $E^1$-term of the spectral sequence is described in terms of the
 homology of $X$. The spectral sequence is different from known spectral
 sequences that compute the homology of iterated loop spaces, such as
 the Eilenberg-Moore spectral sequence and the spectral sequence studied
 by Ahearn and Kuhn in \cite{math/0109041}.  
\end{abstract}

\section{Introduction}
\label{introduction}

The aim of this article is to reveal an unexpected connection
between the combinatorics of the braid arrangements and homological
properties of double (or more highly iterated) loop spaces.

\subsection{Homology of Loop Spaces}

In order to bridge different terminologies, notations, and interests in
these two subjects, hyperplane arrangements and homology of loop spaces,
let us first overview the difficulties in computing the homology of loop
spaces. 

Given a based space $Z$ and a homology theory $h_*(-)$, the homology
$h_*(\Omega Z)$ of the loop space of $Z$ is a functor of $Z$. It would
be nice if we could describe $h_*(\Omega Z)$ in terms of $h_*(Z)$. In
general, we cannot expect such a nice situation, but we do have a
spectral sequence whose $E^2$-term is a functor of $h_*(Z)$. Let us
briefly recall the construction.

Given a diagram of spaces
\begin{equation}
 \begin{diagram}
  \node{}  \node{Y} \arrow{s,r}{g} \\
  \node{X} \arrow{e,t}{f} \node{Z,}
 \end{diagram}
 \label{pullback}
\end{equation}
we can construct a cosimplicial space $\Omega^{\geo}(f,g)$, so-called
the geometric cobar construction as follows. The $k$-th space of
$\Omega^{\geo}(f,g)$ is given by 
\[
 \Omega^{\geo}(f,g)^k = X\times Z^k \times Y
\]
and the structure maps
\begin{eqnarray*}
 d^i & : & \Omega^{\geo}(f,g)^k \longrightarrow \Omega^{\geo}(f,g)^{k+1}
  \\ 
 s^i & : & \Omega^{\geo}(f,g)^{k} \longrightarrow
  \Omega^{\geo}(f,g)^{k-1} 
\end{eqnarray*}
are essentially given by diagonals, projections, $f$, and $g$. The total
space of this cosimplicial space $\Tot(\Omega^{\geo}(f,g))$ is known to
be homeomorphic to the homotopy pullback of the diagram (\ref{pullback})
\[
 \Tot(\Omega^{\geo}(f,g)) \cong \holim(X\rarrow{f} Z \larrow{g} Y).
\]
When $X=Y=\ast$, we obtain a cosimplicial model for the loop space of $Z$
\[
 \Tot(\Omega(\ast,\ast)) \cong \holim(\ast\rarrow{\ast}
 Z\larrow{\ast}\ast) = \Omega Z.
\]

Rector \cite{Rector70} used this cosimplicial model $\Omega^{\geo}(f,g)$
to reformulate the construction of a spectral sequence originally
obtained by Eilenberg and Moore \cite{Eilenberg-Moore66-1,
Eilenberg-Moore66-2}.

When $h_*(-)$ is a multiplicative homology theory satisfying the strong
form of K{\"u}nneth formula 
\[
 h_*(A\times B) \cong h_*(A)\otimes_{h_*} h_*(B),
\]
the $E^2$-term of the spectral sequence for the diagram (\ref{pullback})
can be written as the homology of the algebraic cobar construction
\[
 E^2 \cong H_*(\Omega^{\alg}(f_*,g_*))
\]
where $\Omega^{\alg}(f_*,g_*)$ is a negatively graded chain complex
whose $(-k)$-th term is given by
\[
 \Omega^{\alg}(f_*,g_*)_{-k} = h_*(X)\otimes_{h_*}
 h_*(Z)^{\otimes k} \otimes_{h_*} h_*(Y)
\]
and the boundaries are essentially given by the coalgebra structure on
$h_*(Z)$ and induced maps $f_*$ and $g_*$. 

The homology of the cobar construction $\Omega^{\alg}(f_*,g_*)$ is
isomorphic to the derived functor 
$\Cotor^{h_*(Z)}(h_*(X),h_*(Y))$ of the cotensor product functor over
the coalgebra $h_*(Z)$. Thus we obtain a spectral sequence with 
\[
 E^2 \cong \Cotor^{h_*(Z)}(h_*(X),h_*(Y)).
\]
When $h_*(-)$ is the ordinary homology theory with coefficients in a
field $k$, the spectral sequence is known to converge to the homology of
$\Tot(\Omega^{\geo}(f,g))$ under certain conditions. In particular, we
have a spectral sequence 
\[
 E^2 \cong \Cotor^{H_*(Z;k)}(k,k) \Longrightarrow H_*(\Omega Z;k)
\]
when $Z$ is simply connected.

For other homology theories, especially non-connective ones, however,
the behavior of the spectral sequence could be disastrous and the
$E^{\infty}$-term might not have any relation to the
homology of $\Tot(\Omega^{\geo}(f,g))$ at all.
For example, when $Z$ is the Eilenberg-Mac\,Lane space of
$(\Z/p\Z,2)$-type
\[
 Z=K(\Z/p\Z,2)
\]
and $h_*(-)$ is the mod $p$ $K$-theory
\[
 h_*(-) = K(-;\Z/p\Z),
\]
the $E^2$-term, and hence the $E^{\infty}$-term is trivial but
$h_*(\Omega Z)=h_*(K(\Z/p\Z,1))$ is known to be
nontrivial. The spectral sequence does not give us any information on
$h_*(\Omega Z)$. 

There have been several attempts to find conditions under which the
Eilenberg-Moore spectral sequence can be used to compute the homology of 
$\Omega Z$
\[
 E^2 \cong \Cotor^{h_*(Z)}(h_*,h_*) \Longrightarrow h_*(\Omega Z).
\]
We still don't know a complete answer, but partial answers are
known. One of such answers is given by the author. It is proved in
\cite{Tamaki94} that, when $Z=\Omega^{\ell-1}\Sigma^{\ell} X$, there is a 
spectral sequence with
\begin{equation}
 E^2 \cong \Cotor^{h_*(\Omega^{\ell-1}\Sigma^{\ell} X)}(h_*,h_*),
\label{EMSS} 
\end{equation}
which splits into a direct sum of small spectral sequences each of which
strongly converges to the corresponding summand in the decomposition
\begin{equation}
 \tilde{h}_*(\Omega^{\ell}\Sigma^{\ell} X) \cong \bigoplus_{k}
 \tilde{h}_*(\mathcal{C}_{\ell}(k)_+\wedge_{\Sigma_k} X^{\wedge k}).
\label{homology_decomposition} 
\end{equation}
This is the decomposition induced from the famous stable splitting due
to Snaith \cite{Snaith74}
\begin{equation}
 \Sigma^{\infty}\Omega^{\ell}\Sigma^{\ell} X \simeq \Sigma^{\infty}
 \left(\bigvee_k \mathcal{C}_{\ell}(k)_+\wedge_{\Sigma_k} X^{\wedge
  k}\right).
 \label{Snaith_splitting}
\end{equation}

Thus we do have a spectral sequence that computes $h_*(\Omega Z)$ in
this case. It is also proved by the author \cite{Tamaki06} that the
spectral sequence (\ref{EMSS}) is isomorphic to the Eilenberg-Moore
spectral sequence. 

The space $\mathcal{C}_{\ell}(k)$ appeared in the direct sum
decomposition (\ref{homology_decomposition}) is
the space of little cubes studied by May in \cite{May72}.
The author constructed the spectral sequence (\ref{EMSS}) by defining a
filtration on each $\mathcal{C}_{\ell}(k)$
\begin{equation}
 \emptyset = F_{-k-1}\mathcal{C}_{\ell}(k) \subset
 F_{-k}\mathcal{C}_{\ell}(k) \subset \cdots \subset
 F_{-2}\mathcal{C}_{\ell}(k) \subset F_{-1}\mathcal{C}_{\ell}(k) =
 F_0\mathcal{C}_{\ell}(k) = \mathcal{C}_{\ell}(k),
 \label{gravity_filtration}
\end{equation}
which is a completely different method from the cosimplicial
construction of the Eilenberg-Moore spectral sequence.

\subsection{The Little $2$-Cubes and the Braid Arrangements}

The space $\mathcal{C}_{\ell}(k)$ of $k$ little $\ell$-cubes is homotopy
equivalent to the configuration space $F(\R^{\ell},k)$ of distinct $k$
points in $\R^{\ell}$. When $\ell=2$
\[
 \mathcal{C}_2(k) \simeq F(\C,k) =
 \lset{(z_1,\cdots,z_k)\in\C^k}{z_i\neq z_{j} \text{ if } i\neq j}.
\]
For $1\le i < j\le k $, define a hyperplane in $\R^k$
\[
 L_{i,j} = \lset{(x_1,\cdots, x_k)\in \R^k}{x_i=x_j}.
\]
Then $F(\C,k)$ is the complement of the complexification of the
real central hyperplane arrangement $\lset{L_{i,j}}{1\le i<j\le k}$
\[
 F(\C,k) = \C^k - \bigcup_{1\le i<j\le k} L_{i,j}\otimes\C.
\]
Hyperplanes in this arrangement contains a line $\{x_1=\cdots=x_k\}$ in
common. To make it essential, let 
\[
 \mathfrak{h}_k = \lset{(x_1,\cdots,x_k)\in\R^{k}}{x_1+\cdots+x_k=0}
\]
and 
\[
 L'_{i,j} = L_{i,j}\cap \mathfrak{h}_k.
\]
The arrangement $\mathcal{A}_{k-1} = \{L'_{i,j} \mid 1\le i<j\le k\}$ is
essential and we have homotopy equivalences
\begin{eqnarray*}
 \mathcal{C}_2(k) & \simeq & F(\C,k) \\
 & = & \C^k - \bigcup_{1\le i<j\le k} L_{i,j}\otimes\C \\
 & \simeq & \mathfrak{h}_k\otimes\C -\bigcup_{1\le i<j\le k}
  L'_{i,j}\otimes \C. 
\end{eqnarray*}
Note that these homotopy equivalences respect the right action of the
symmetric group $\Sigma_k$.

The homotopy types of complements of complexified arrangements have been
actively studied by many people. One of the most interesting and useful
constructions in this subject is a result of M.~Salvetti
\cite{Salvetti87}, who found a construction of a finite simplicial
complex $\Sal(\mathcal{A})$ for a real central essential arrangement  
$\mathcal{A}$, embedded in the complement
of the complexification of $\mathcal{A}$ as a deformation retract. Thus
we obtain
\begin{equation}
 \mathcal{C}_2(k) \simeq \mathfrak{h}_k\otimes\C -\bigcup_{1\le i<j\le k}
  L'_{i,j}\otimes \C \simeq \Sal(\mathcal{A}_{k-1}).
  \label{C_and_Sal}
\end{equation}

There should be a combinatorial meaning of the filtration
(\ref{gravity_filtration}) on
$\mathcal{C}_2(k)$ when translated to the Salvetti 
complex $\Sal(\mathcal{A}_{k-1})$ of the braid arrangement under the
homotopy equivalence (\ref{C_and_Sal}).

\subsection{Statement of Results}

For any real central hyperplane arrangement $\mathcal{A}$, Salvetti, on
the other hand, defined a structure of regular cell complex on
$\Sal(\mathcal{A})$ by combining  simplices by using the combinatorial
structure of the face lattice of $\mathcal{A}$.

It turns out that the filtration (\ref{gravity_filtration}) on
$\mathcal{C}_2(k)$ designed for the Eilenberg-Moore spectral sequence 
(\ref{EMSS}) coincides with the filtration defined by the cellular
structure on $\Sal(\mathcal{A}_{k-1})$ under the homotopy equivalence 
(\ref{C_and_Sal}). 

\begin{introtheorem}[Theorem \ref{Subquotient}]
 Let $F_{-s}\Sal(\mathcal{A}_{k-1})$ be the $(k-s)$-skeleton of
 $\Sal(\mathcal{A}_{k-1})$ under the cellular structure defined in
 \cite{Salvetti87}. Then the homotopy equivalence
 \[
  \varphi_k : \Sal(\mathcal{A}_{k-1}) \rarrow{\simeq} \mathcal{C}_{2}(k) 
 \]
 in (\ref{C_and_Sal}) preserves filtrations
 and induces a $\Sigma_k$-equivariant homotopy
 equivalence on each subquotient
 \[
 \varphi_k :
 F_{-s}\Sal(\mathcal{A}_{k-1})/F_{-s-1}\Sal(\mathcal{A}_{k-1})
 \rarrow{\simeq} F_{-s}\mathcal{C}_2(k)/F_{-s-1}\mathcal{C}_2(k).
 \]
\end{introtheorem}

A filtration on each $\mathcal{C}_2(k)$ induces a spectral sequence
converging to $h_*(\mathcal{C}_2(k)_{+}\wedge_{\Sigma_k} X^{\wedge k})$
and the Eilenberg-Moore spectral sequence (\ref{EMSS}) decomposes into a
direct sum of these small spectral sequences. Let
$\{E^r(\mathcal{A}_{k-1})\}$ be the spectral sequence for
$h_*(\Sal(\mathcal{A}_{k-1})_{+}\wedge_{\Sigma_k} X^{\wedge k})$ defined
by the skeletal filtration on $\Sal(\mathcal{A}_{k-1})$. Then the
$E^1$-term of the spectral sequence defined in \cite{Tamaki94} for
$h_*(\Omega^2\Sigma^2 X)$ can be identified as follows
\begin{eqnarray*}
 E^1_{-s,t} & = & \bigoplus_{k\ge 0} E^1_{-s,t}(\mathcal{A}_{k-1}) \\
 & \cong &
  \bigoplus_{k\ge 0}C_{k-s}(\Sal(\mathcal{A}_{k-1})) \otimes_{\Sigma_k}
  \tilde{h}_{t-k}\left(X^{\wedge k}\right) \\
 & \cong & \bigoplus_{k\ge 0}C_{k-s}(\Sal(\mathcal{A}_{k-1}))
  \otimes_{\Sigma_k} \tilde{h}_{t}\left((\Sigma X)^{\wedge k}\right),
\end{eqnarray*}
where $C_*(-)$ denotes the cellular chain complex functor. This is an
isomorphism of chain complexes. When $h_*(-)$
is multiplicative and satisfies the strong form of K\"unneth formula, we
have 
\[
 E^1_{-s,*} \cong \bigoplus_{k\ge 0}
 C_{k-s}(\Sal(\mathcal{A}_{k-1}))\otimes_{\Sigma_k} \tilde{h}_*(\Sigma
 X)^{\otimes k}.
\]
The $(k-s)$-cells in $\Sal(\mathcal{A}_{k-1})$ are labeled by a pair of 
an ordered partition $\lambda$ and a permutation $\sigma$ which is a
subdivision of $\lambda$, in which case we denote $\lambda \le
\sigma$. Then we have 
\begin{eqnarray*}
 E^1_{-s,*} & \cong & \bigoplus_{k\ge 0}
 \Z\langle [D(\lambda,\sigma)] \mid \lambda \in
   \Pi_{k,k-s}, \sigma\in\Sigma_k, \lambda\le
 \sigma\rangle\otimes_{\Sigma_k} \tilde{h}_*(\Sigma X)^{\otimes k} \\
 & \cong & \bigoplus_{k\ge 0} \Z\langle [D(\lambda, (1|\cdots|k))] \mid
  \lambda \in O_{k,k-s}\rangle\otimes \tilde{h}_*(\Sigma X)^{\otimes k},
\end{eqnarray*}
where $(1|\cdots|k) \in \Sigma_k$ is the identity and $\Pi_{k,k-s}$ and
$O_{k,k-s}$ are the set of ordered partitions and order-preserving
partitions of rank $k-s$, respectively. See \S\ref{Configuration} for
precise definitions.

On the other hand, it is proved in \cite{Tamaki06} that the $E^1$-term
of the spectral sequence is isomorphic to the cobar complex of
$h_*(\Omega\Sigma^2 X) \cong T(\tilde{h}_*(\Sigma X))$, the tensor
algebra on $\tilde{h}_*(\Sigma X)$. This means that the cellular
structure of the Salvetti complex describes the iterated cobar
construction. This immediately gives us the description of the cobar
differential stated in Introduction of \cite{Tamaki06}. (For the precise
definitions, see p.\pageref{OrderedPartition} and p.\pageref{d^1}.)

\begin{introtheorem}[Corollary \ref{d^1}]
 The first differential $d^1$ on the spectral sequence is given by the
 following formula:
 \begin{multline*}
  d^1_{-s,*}([D(\lambda,(1|\cdots|k))]\otimes[x_1|\cdots|x_k]) \\
  = \sum_{\tau \in O_{k,k-s-1}, \lambda<\tau} [D(\tau,(1|\cdots|k))]\otimes
 \left(\sum_{\rho\in S_{t(\tau)}}
 \sgn(\rho)[x_1|\cdots|x_k]\cdot\rho\right),
 \end{multline*}
 where $S_{t(\tau)}$ is the set of shuffles of the same type as $\tau$.
\end{introtheorem}

With this description, it is easy to see the first differentials are
induced by space-level shuffles.
\[
 \bigvee_{\lambda\in O_{k,k-s}} S^{k-s}_{(\lambda,(1|\cdots|k))} \wedge
 (\Sigma X)^{\wedge k} 
 \longrightarrow \Sigma\left(\bigvee_{\tau\in O_{k,k-s-1, \lambda<\tau}}
 S^{k-s-1}_{(\tau,(1|\cdots|k))}
 \wedge (\Sigma X)^{\wedge k}\right),
\]
where $S^{k-s}_{(\lambda,(1|\cdots|k))}$ and
$S^{k-s-1}_{(\tau,(1|\cdots|k))}$ are copies of spheres $S^{k-s}$ and
$S^{k-s-1}$, respectively.

This map may be of some use to study maps between wedge powers of
suspended spaces.

The spectral sequence (\ref{EMSS}) is defined for not only $\ell=2$ but
all $\ell\ge 1$. On the other hand, it was observed by Bj{\"o}rner and
Ziegler \cite{Bjorner-Ziegler92} that the construction of
$\Sal(\mathcal{A})$ can be extended to subspace arrangements of higher
codimensions by using the notion of oriented $k$-matroid. For a real
central essential arrangement $\mathcal{A}$ in a vector space $V$, there
is a simplicial complex $\Sal^{(\ell)}(\mathcal{A})$ embedded in the
complement of the $\ell$-dimensionalization
\[
 \Sal^{(\ell)}(\mathcal{A}) \hookrightarrow V\otimes\R^{\ell} -
 \bigcup_{L\in\mathcal{A}} L\otimes\R^{\ell} 
\]
as a deformation retract. We can generalize the cellular structure of
$\Sal(\mathcal{A}) = \Sal^{(2)}(\mathcal{A})$ to
$\Sal^{(\ell)}(\mathcal{A})$ to make it a regular cell complex. Under
the stable homotopy equivalence
\[
 \Omega^{\ell}\Sigma^{\ell} X \shot \bigvee_{k}
 \mathcal{C}_{\ell}(k)_{+}\wedge X^{\wedge k} \simeq \bigvee_k
 F(\R^{\ell},k)_{+}\wedge X^{\wedge k} \simeq 
 \bigvee_k \Sal^{(\ell)}(\mathcal{A}_{k-1})_{+}\wedge X^{\wedge k}, 
\]
the skeletal filtration on $\Sal^{(\ell)}(\mathcal{A}_{k-1})$ induces a
spectral sequence that computes $h_*(\Omega^{\ell}\Sigma^{\ell}X)$.

\begin{introtheorem}[Theorem \ref{higher_gravity_spectral_sequence}]
 There exists a spectral sequence for any homology theory
 \begin{equation}
  E^1 \cong \bigoplus_k
 C_*(\Sal^{(\ell)}(\mathcal{A}_{k-1}))\otimes_{\Sigma_k} 
 \tilde{h}_*(X^{\wedge k}) \Longrightarrow
 h_*(\Omega^{\ell}\Sigma^{\ell} X),
 \label{higher_GSS}
 \end{equation}
 which is a direct sum of spectral sequences each of which strongly
 converges to the corresponding direct summand in
 (\ref{homology_decomposition}). 
\end{introtheorem}

The spectral sequence is different from the Eilenberg-Moore spectral
sequence (\ref{EMSS}) when $\ell>2$. It is also different from the
spectral sequence studied by Ahearn and Kuhn \cite{math/0109041}. Thus
we obtain a new spectral sequence which can be used to compute
$h_*(\Omega^{\ell}\Sigma^{\ell} X)$. As we will see in \S\ref{remarks},
our spectral sequence is much finer than their spectral sequence.

\bigskip

The organization of this paper is as follows: we recall
the construction of the Salvetti complex in \S\ref{Complex}.
Salvetti's construction has been intensively studied and
several alternative descriptions and interpretations are known. We use
the description in terms of matroid product developed in
\cite{Gelfand-Rybnikov89}, \cite{Arvola91}, and
\cite{Bjorner-Ziegler92}. The properties of the 
filtration on the space of little
cubes are summarized in \S\ref{LittleCubes}. After these
preliminary sections, we compare the Salvetti complex for the braid
arrangement $\mathcal{A}_{k-1}$ and the filtration on the space of $k$
little $2$-cubes in \S\ref{Configuration} and prove Theorem
\ref{Subquotient}. Higher order analogues are discussed in 
\S\ref{higher_Salvetti}, where we construct the spectral sequence
(\ref{higher_GSS}) by using oriented matroids. Comparisons with other
spectral sequences that compute the homology of iterated loop spaces are
given in \S\ref{remarks}.

\begin{acknowledgement}
 It was an e-mail message from Ryan Budney
 which drove the author to try to find a relation between the Salvetti
 complex and the filtration on little $2$-cubes. The author is grateful
 to him for his intriguing question. The result of this paper
 partially answers his question. Tilman Bauer kindly informed me of the
 inaccuracy of the statement of the main result in
 \cite{Tamaki94}. Finally the author would like to thank Fred 
 Cohen and Jim Stasheff for their interests and suggestions.
\end{acknowledgement}

\section{The Salvetti Complex}
\label{Complex}

For a real central arrangement $\mathcal{A}$ in $\R^n$, Salvetti defined a
simplicial complex $\Sal(\mathcal{A})$ embedded in the complement of the
complexification of the arrangement in $\C^n$. Let us recall the
construction of $\Sal(\mathcal{A})$.  Although the combinatorial
informations on a real central arrangement can be translated into the
language of oriented matroids and the construction of the Salvetti
complex can be generalized to oriented matroids, we follow the original
treatment \cite{Salvetti87} in this section, since we need explicit
descriptions of vertices for later use. We use the matroid description
of the Salvetti complex in \S\ref{higher_Salvetti}.

Let us first fix notations and terminologies. Let $\mathcal{A} =
\{M_j\}_{j\in J}$ be a
central arrangement of hyperplanes in $\R^n$. The arrangement
$\mathcal{A}$ defines a stratification of $\R^n$.
\begin{eqnarray*}
 M^0 & = & \R^n-\bigcup_{j\in J} M_j \\
 M^1 & = & \bigcup_{j\in J} \left(M_j- \bigcup_{k\neq j} M_j\cap
			     M_k\right) \\
 & \vdots & \\
 M^{|J|} & = & \bigcap_{j\in J} M_j \\
 \R^n & = & M^0\cup M^1\cup \cdots \cup M^{|J|}
\end{eqnarray*}

The connected components of each stratum are called faces and a face
in the top stratum is called a chamber. The sets of all faces and of
faces of codimension $i$ are denoted by
\begin{eqnarray*}
 \mathcal{F}(\mathcal{A}) & = & \{\text{faces of }\mathcal{A}\} \\
 \mathcal{F}^i(\mathcal{A}) & = & \{F\in\mathcal{F}(\mathcal{A}) \mid
  \codim F = i\}.
\end{eqnarray*}
In particular, $\mathcal{F}^0(\mathcal{A})$ is the set of chambers.
We need an ordering in $\mathcal{F}(\mathcal{A})$. Although a popular
way to define an ordering is by the ``reverse inclusion'', Salvetti
defines as follows
\[
 F\ge F' \Longleftrightarrow \overline{F} \supset F'.
\]

In order to define a simplicial complex, we need vertices.
For each face $F$, we choose a point $w(F) \in F$ and the collections of
these points are denoted by
\begin{eqnarray*}
 \mathcal{V}(\mathcal{A}) & = & \{w(F) \mid
  F\in\mathcal{F}(\mathcal{A})\} \\
 \mathcal{V}^i(\mathcal{A}) & = & \{w(F) \mid
  F\in\mathcal{F}^i(\mathcal{A})\}.
\end{eqnarray*}
These are points in $\R^n$. We need to add some imaginary coordinates
to obtain points in $\C^n$. Salvetti defines the imaginary part by using
the following fact.

\begin{lemma}
 For $v\in \mathcal{V}^0(\mathcal{A})$ and
 $F \in \mathcal{F}(\mathcal{A})$, there exists a unique point
 $w(v,F) \in \mathcal{V}^0(\mathcal{A})$ with the following properties:
 \begin{enumerate}
  \item $v$ and $w(v,F)$ belong to the same chamber of
	$\mathcal{A}_{\supset F}$, which is an arrangement defined by
	\[
	\mathcal{A}_{\supset F} = \{H \in \mathcal{A} \mid H\supset F\}.
	\]
  \item $w(v,F)$ belongs to a chamber in
	$\{C\in\mathcal{F}^0(\mathcal{A})\mid C\ge F\}$.
 \end{enumerate}
\end{lemma}

\begin{figure}[h]
\begin{center}
 \begin{picture}(170,160)(-20,-10)
  \put(0,0){\line(1,1){130}}
  \put(133,130){\makebox(0,0)[bl]{$H_1$}}

  \put(130,0){\line(-1,1){130}}
  \put(-3,130){\makebox(0,0)[br]{$H_2$}}

  \put(115,65){\circle*{3}}
  \put(118,65){\makebox(0,0)[l]{$v$}}

  \put(58,30){\circle*{3}}
  \put(61,30){\makebox(0,0)[l]{$w(v,F)$}}
  \put(58,10){\makebox(0,0){$C$}}

  \thicklines
  \put(0,0){\line(1,1){65}}
  \put(30,40){\makebox(0,0)[br]{$F$}}

  \thinlines
  \put(12,0){\line(0,1){12}}
  \put(24,0){\line(0,1){24}}
  \put(36,0){\line(0,1){36}}
  \put(48,0){\line(0,1){48}}
  \put(60,0){\line(0,1){60}}
  \put(72,0){\line(0,1){72}}
  \put(84,0){\line(0,1){84}}
  \put(96,0){\line(0,1){96}}
  \put(108,0){\line(0,1){108}}
  \put(120,0){\line(0,1){120}}
 \end{picture}
\end{center}
 \caption{$w(v,F)$}
\end{figure}

The vertices of the Salvetti complex are given by
\[
 w(F) + i(w(v,F)-w(F))
\]
for each element $v \in \mathcal{V}^0(\mathcal{A})$ and a face $F$.
Note that $w(v,F)$ is a vertex corresponding to a chamber. Thus the
vertices of the Salvetti complex are given by
\[
 v(F,C) = w(F)+i(w(C)-w(F))
\]
for a face $F$ and a chamber $C$ satisfying a certain condition. When
$v$ varies, all points in the chamber $C$ with $C\ge F$ appear as
$w(v,F)$. Thus we have the following.

\begin{definition}
 For a real central arrangement $\mathcal{A}$, define the vertex set by
 \[
 \sk_0(\Sal(\mathcal{A})) = \{v(F,C) \mid F \in
 \mathcal{F}(\mathcal{A}), C\in \mathcal{F}^0(\mathcal{A}), F\le C\}
 \]
\end{definition}

In oder to define higher dimensional simplices, we need to specify a
condition for a collection of vertices to form a simplex.  The condition
can be described by using a ``face-chamber''pairing
\[
  \mathcal{F}(\mathcal{A})\times\mathcal{F}^0(\mathcal{A})
  \longrightarrow \mathcal{F}^0(\mathcal{A}).
\]
This can be extended to a ``face-face'' pairing
\[
  \mathcal{F}(\mathcal{A})\times\mathcal{F}(\mathcal{A})
  \longrightarrow \mathcal{F}(\mathcal{A}).
\]

\begin{proposition}
 \label{matroid_product_of_faces}
 For $F, G \in \mathcal{F}(\mathcal{A})$, there exists a unique face
 $F\circ G \in  \mathcal{F}(\mathcal{A})$ with the following
 properties: let $H \in \mathcal{F}(\mathcal{A}_{\supset |F|})$ be the
 unique face with $G\subset H$. Then $F\circ G$ is the unique face
 with $F\circ G \subset H$ and $F\le F\circ G$.
\end{proposition}

\begin{lemma}[\cite{Arvola91}]
 \label{PairingProperties}
 The above face-face pairing has the following properties:
 \begin{enumerate}
  \item The pairing is associative.
  \item If $G$ is a chamber, $F\circ G$ is also a chamber.
  \item If $F$ is a chamber, $F\circ G = F$.
  \item If $G\le G'$, then $F\circ G\le F\circ G'$.
 \end{enumerate}
\end{lemma}

With this pairing, the simplices of the Salvetti complex can be defined
as follows.

\begin{definition}
 For a chain $F_0>F_1>\cdots>F_j$ in $\mathcal{F}(\mathcal{A})$ and a
 chamber $C$ with $C \ge F_j$, define a simplex by
 \[
 s(F_0,\cdots, F_j;C) = \Conv(\{v(F_0,F_0\circ C), \cdots,
 v(F_j,F_j\circ C)\})
 \]
 where $\Conv(S)$ denotes the convex hull of $S$.

 The Salvetti complex is defined by
 \[
  \Sal(\mathcal{A}) = \bigcup_{j}\{s(F_0,\cdots, F_j;C) \mid F_0>\cdots>F_j,
   C\ge F_j, C\in\mathcal{F}^0(\mathcal{A})\}.
 \]
\end{definition}

Salvetti proves that $\Sal(\mathcal{A})$ is a deformation retract of the
complement of the complexification $\mathcal{A}^{\C}$ of $\mathcal{A}$.

\begin{theorem}[Salvetti \cite{Salvetti87}]
 $\Sal(\mathcal{A})$ is contained in the complement of the
 complexification of $\mathcal{A}$ and the inclusion is a homotopy
 equivalence
 \[
  \Sal(\mathcal{A}) \simeq \C^n -\bigcup_{M\in \mathcal{A}} M^{\C}
 \]
 where $M^{\C}$ is the complexification of $M$.
\end{theorem}

Although $\Sal(\mathcal{A})$ has a nice simplicial structure, it is more
convenient to combine simplices which form a cell. In fact, Salvetti
defines a CW-structure. Using the facet-facet pairing, we have the
following description.

\begin{definition}
 For $F \in \mathcal{F}(\mathcal{A})$ and $C \in
 \mathcal{F}^0(\mathcal{A})$ with $F\le C$, define a subset of
 $\sk_0(\Sal(\mathcal{A}))$ by
 \[
  \mathcal{D}(F,C) = \{v(G,G\circ C) \mid G\ge F\}.
 \]
 This set is regarded as a poset by
 \[
  v(G,G\circ C) \le v(H,H\circ C) \Longleftrightarrow G\le H.
 \]

 The (geometric realization of the) order complex of $\mathcal{D}(F,C)$
 is denoted by $D(F,C)$. 
\end{definition}

Salvetti proves the following. See \cite{Arvola91}.

\begin{lemma}
 \label{D(F,C)}
 The complex $D(F,C)$ has the following properties:
 \begin{enumerate}
  \item The inclusion of vertices induces a simplicial embedding
	\[
	 D(F,C) \hookrightarrow \Sal(\mathcal{A}).
	\]
  \item $D(F,C)$ is homeomorphic to a disk of dimension $\codim F$.
  \item The boundary of $D(F,C)$ is given by
	\[
	 \partial D(F,C) = \bigcup_{G>F} D(G,G\circ C).
	\]
  \item The decomposition
	\[
	 \Sal(\mathcal{A}) = \bigcup_{v(F,C)\in
	\sk_0(\Sal(\mathcal{A}))} (D(F,C)-\partial D(F,C))
	\]
	gives a structure of a finite regular cell complex.
 \end{enumerate}
\end{lemma}
\section{The Gravity Filtration on Little Cubes}
\label{LittleCubes}

In \cite{Tamaki94}, the author introduced a filtration on the space of
$k$ little $\ell$-cubes, $\mathcal{C}_{\ell}(k)$, for all $\ell$ and
$k$. The subquotients
$F_{-s}\mathcal{C}_{\ell}(k)/F_{-s-1}\mathcal{C}_{\ell}(k)$ are 
analyzed in \cite{Tamaki06}. In this section, we recall the properties
of this filtration.

$\mathcal{C}_{\ell}(k)$ is the space of little cubes in $I^{\ell}$. In
order to compare it with the Salvetti complex in the next section,
however, it is more convenient to consider cubes in $\R^{\ell}$.

\begin{convention}
 In the rest of this paper, $\mathcal{C}_{\ell}(k)$ denotes the space of
 $k$ little $\ell$-cubes in $\R^{\ell}$ whose images have disjoint
 interiors to each other and whose edges are parallel to the coordinate 
 axes. However, when we draw a picture, we use the usual picture 
 for little cubes, i.e.\ ``non-overlapping small boxes in a big box''.
\end{convention}

In order to state the result in \cite{Tamaki06}, we need the following
notations.
\begin{definition}
 Fix $\varepsilon >0$, and define an $(k-1)$-dimensional convex polytope
 $P^k$ in the hyperplane
 $\mathfrak{h}_k = \{(t_1,\cdots, t_k)\in\R^{k} \mid t_1+\cdots+t_k = 0\}$
 by 
 \[
  P^{k} = \{(t_1,\cdots, t_k)\in\mathfrak{h}_k \mid
 |t_i-t_j|<\varepsilon \text{ for } i\neq j\}.
 \]
 For a subset $S=\{s_1,\cdots,s_{\ell}\} \subset \{1,\cdots, k\}$, define
 \begin{eqnarray*}
  \R^S & = & \Map(S,\R) = \{(t_{s_1},\cdots, t_{s_{\ell}}) \mid
   t_{s_1},\cdots,t_{s_{\ell}}\in \R\} \\
  \mathfrak{h}_S & = & \{(t_{s_1},\cdots,t_{s_{\ell}})\in \R^S \mid
   t_{s_1}+\cdots+t_{s_{\ell}} = 0\}
 \end{eqnarray*}
 and define
 \[
  P^{S} = \{(t_{s_1},\cdots, t_{s_{\ell}})\in\mathfrak{h}_S \mid
 |t_{s}-t_{s'}|<\varepsilon \text{ for } s\neq s'\}.
 \]
 Also for little cubes, we denote
 \[
  \mathcal{C}_{\ell}(S) = \{(c_{s_1},\cdots,c_{s_{\ell}}) \in
 \Map(S,\mathcal{C}_{\ell}(1)) \mid c_{s_i}(\Int I^{\ell}) \cap
 c_{s_j}(\Int I^{\ell}) = \emptyset\}.
 \]
\end{definition}

We also need the notations for partitions. Here a partition always means
an ordered partition.

\begin{definition}
 A partition of $\{1,\cdots, k\}$ is a surjective map
 \[
 \lambda : \{1,\cdots,k\} \longrightarrow \{1,\cdots,k-r\}
 \]
 for some $0\le r < k$. The number $r$ is called the rank of this
 partition.
 
 The set of partitions of $\{1,\cdots,k\}$ is denoted by
 $\Pi_k$. The subset of rank $r$ partitions is denoted by
 $\Pi_{k,r}$.

 Note that rank $0$ partitions are nothing but elements of $\Sigma_k$.
\end{definition}

With these notations, it is proved in \cite{Tamaki06} that

\begin{proposition}
 \label{SubquotientsForC}
 We have the following $\Sigma_k$-equivariant homotopy equivalence
 \begin{multline*}
  F_{-s}\mathcal{C}_{\ell}(k)/F_{-s-1}\mathcal{C}_{\ell}(k) \simeq \\
 \bigvee_{\lambda \in \Pi_{k,k-s}}
  \mathcal{C}_{\ell-1}({\lambda^{-1}(1)})_{+} \wedge
 (P^{\lambda^{-1}(1)}/\partial P^{\lambda^{-1}(1)}) \wedge \cdots \wedge 
 \mathcal{C}_{\ell-1}({\lambda^{-1}(s)})_{+} \wedge
  (P^{\lambda^{-1}(s)}/\partial P^{\lambda^{-1}(s)}) \\
  = \bigvee_{\lambda \in \Pi_{k,k-s}}
  \mathcal{C}_{\ell-1}(\lambda)_+ \wedge P^{\lambda}/\partial P^{\lambda},
 \end{multline*}
 where
 \begin{eqnarray*}
  \mathcal{C}_{\ell-1}(\lambda) & = &
   \mathcal{C}_{\ell-1}({\lambda^{-1}(1))}\times \cdots\times
   \mathcal{C}_{\ell-1}({\lambda^{-1}(s)}) \\
  P^{\lambda} & = & P^{\lambda^{-1}(1)}\times \cdots
  \times P^{\lambda^{-1}(s)}
 \end{eqnarray*}
 for a partition $\lambda$ of rank $k-s$.
\end{proposition}

This homotopy equivalence is proved as follows. The first step is to
replace $F_{-s}\mathcal{C}_{\ell}(k)$ by horizontally decomposable cubes.

\begin{definition}
 Let $\mathcal{D}_{\ell}^s(k)$ be the subset of $\mathcal{C}_{\ell}(k)$
 consisting of cubes which are horizontally decomposable into $s$
 collections. 

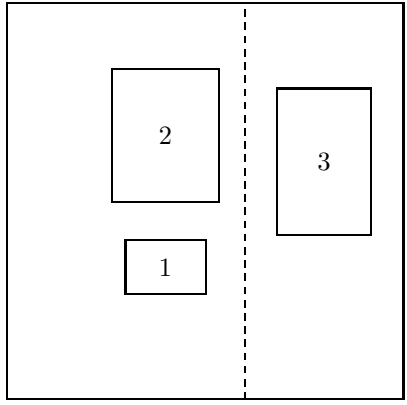
\begin{center}
 \begin{picture}(150,150)(0,0)
  \put(0,0){\line(1,0){150}}
  \put(0,0){\line(0,1){150}}
  \put(150,0){\line(0,1){150}}
  \put(0,150){\line(1,0){150}}

  \put(60,50){\makebox(0,0){\framebox(30,20){$1$}}}
  \put(60,100){\makebox(0,0){\framebox(40,50){$2$}}}
  \put(120,90){\makebox(0,0){\framebox(35,55){$3$}}}

  \multiput(90,0)(0,5){30}{\line(0,1){2.5}}
 \end{picture}
\end{center}
\end{definition}

More precise definition can be given in terms of the operad structure
map. See \cite{Tamaki06}, for details.

Note that $\mathcal{D}_{\ell}^s(k)$ is included in
$F_{-s}\mathcal{C}_{\ell}(k)$. 

\begin{lemma}
 \label{FirstDeformation}
 The inclusion induces a $\Sigma_k$-equivariant homotopy equivalence
 \[
 \mathcal{D}_{\ell}^s(k)/\mathcal{D}_{\ell}^{s+1}(k) \simeq 
 F_{-s}\mathcal{C}_{\ell}(k)/F_{-s-1}\mathcal{C}_{\ell}(k).
 \]

 We also have
 \begin{equation}
 \mathcal{D}_{\ell}^s(k)/\mathcal{D}_{\ell}^{s+1}(k) \simeq
 \bigvee_{\lambda\in \Pi_{k,s}}
 \mathcal{D}_{\ell}^1(\lambda^{-1}(1))/
 \mathcal{D}_{\ell}^{2}(\lambda^{-1}(1))  
 \wedge \cdots\wedge
 \mathcal{D}_{\ell}^1(\lambda^{-1}(s))/
 \mathcal{D}_{\ell}^{2}(\lambda^{-1}(s)).  
 \label{DecompositionOfD}
 \end{equation}
\end{lemma}

Thus it is enough to analyze
$\mathcal{D}_{\ell}^1(k)/\mathcal{D}_{\ell}^2(k)$. We 
can make $\mathcal{D}_{\ell}^1(k)$ even smaller by removing unnecessary
cubes. Note that the following element in $\mathcal{D}^1(3)$ can be
moved into $\mathcal{D}^2(3)$ by shrinking the first coordinates of
cubes.
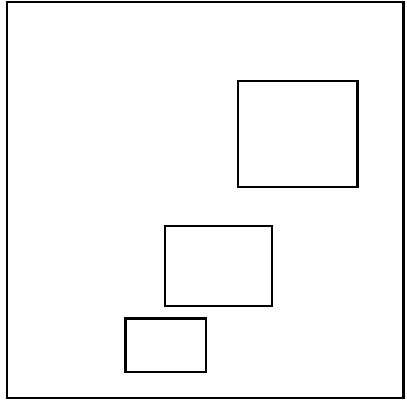
\begin{center}
 \begin{picture}(150,150)(0,0)
  \put(0,0){\line(1,0){150}}
  \put(0,0){\line(0,1){150}}
  \put(150,0){\line(0,1){150}}
  \put(0,150){\line(1,0){150}}

  \put(60,20){\makebox(0,0){\framebox(30,20){}}}
  \put(80,50){\makebox(0,0){\framebox(40,30){}}}
  \put(110,100){\makebox(0,0){\framebox(45,40){}}}
 \end{picture}
\end{center}

All we need is cubes in $\mathcal{D}_{\ell}^1(k)$ which can be skewered
vertically.

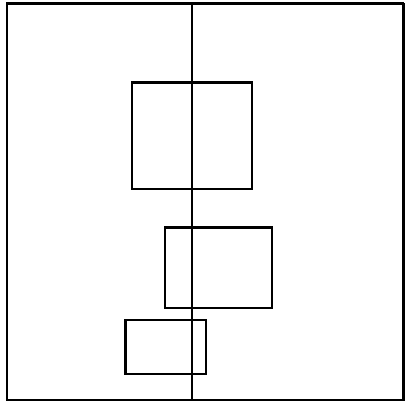
\begin{center}
 \begin{picture}(150,150)(0,0)
  \put(0,0){\line(1,0){150}}
  \put(0,0){\line(0,1){150}}
  \put(150,0){\line(0,1){150}}
  \put(0,150){\line(1,0){150}}

  \put(60,20){\makebox(0,0){\framebox(30,20){}}}
  \put(80,50){\makebox(0,0){\framebox(40,30){}}}
  \put(70,100){\makebox(0,0){\framebox(45,40){}}}

  \put(70,0){\line(0,1){150}}
 \end{picture}
\end{center}

\begin{definition}
 Define $G_{-s}\mathcal{C}_{\ell}(k)$ to be the subset of
 $\mathcal{D}_{\ell}^s(k)$ 
 consisting of cubes $(c_1,\cdots, c_j)$ which cannot be decomposed into
 $(s-1)$ collections of cubes each of which can be skewered by a
 vertical line (hyperplane) intersecting with each interior.
\end{definition}

\begin{lemma}
 The inclusion induces a $\Sigma_k$-equivariant homotopy equivalence
 \[
  G_{-1}\mathcal{C}_{\ell}(k)/G_{-2}\mathcal{C}_{\ell}(k) \simeq
 \mathcal{D}_{\ell}^1(k)/\mathcal{D}_{\ell}^2(k).
 \]
\end{lemma}

Finally we adjust the radii of the first coordinates of cubes in
$G_{-1}\mathcal{C}_{\ell}(k)$ so that they have the same fixed radii.

\begin{definition}
 For $\varepsilon>0$, let $\mathcal{C}_{\ell}^{\varepsilon}(k)$ be the
 subspace 
 of $\mathcal{C}_{\ell}(k)$ consisting of cubes having the radii of the
 first coordinates $\varepsilon$.
\end{definition}

\begin{lemma}
 The inclusion induces a $\Sigma_k$-equivariant homotopy equivalence
 \[
 G_{-1}\mathcal{C}_{\ell}(k)\cap\mathcal{C}_{\ell}^{\varepsilon}(k)/
 G_{-2}\mathcal{C}_{\ell}(k)\cap\mathcal{C}_{\ell}^{\varepsilon}(k)
 \simeq G_{-1}\mathcal{C}_{\ell}(k)/G_{-2}\mathcal{C}_{\ell}(k).
 \]
\end{lemma}

Now it is easy to describe
$G_{-1}\mathcal{C}_{\ell}(k)\cap\mathcal{C}_{\ell}^{\varepsilon}(k)/
 G_{-2}\mathcal{C}_{\ell}(k)\cap\mathcal{C}_{\ell}^{\varepsilon}(k)$.

\begin{lemma}
 Define
 \begin{eqnarray*}
  \widetilde{P}^{k} & = & \{(t_1,\cdots, t_k)\in\R^k \mid
   |t_i-t_j|<\varepsilon \text{ for all } i, j\} \\
  d\widetilde{P}^k & = & \{(t_1,\cdots, t_k)\in\R^k \mid
   |t_i-t_j|=\varepsilon \text{ for some } i\neq j\}.
 \end{eqnarray*}
 Then we have a $\Sigma_k$-equivariant homeomorphism
 \[
 G_{-1}\mathcal{C}_{\ell}(k)\cap\mathcal{C}_{\ell}^{\varepsilon}(k)/
 G_{-2}\mathcal{C}_{\ell}(k)\cap\mathcal{C}_{\ell}^{\varepsilon}(k) \cong
 \widetilde{P}^k/d\widetilde{P}^k\wedge \mathcal{C}_{\ell-1}(k)_{+}. 
 \]
\end{lemma}

Under the decomposition of the permutation representation
\begin{equation}
 \R^k \cong \mathfrak{h}_k\oplus \{(t,\cdots,t)\mid t\in\R\},
  \label{PermutationRepresentation}
\end{equation}
we see that the inclusion induces a $\Sigma_k$-equivariant homotopy
equivalence
\[
 P^k/\partial P^k \simeq \widetilde{P}^k/d\widetilde{P}^k.
\]

The homotopy equivalence in Proposition \ref{SubquotientsForC} is
induced by the composition of the inclusions
\begin{equation}
 P^{\lambda^{-1}(i)} \hookrightarrow \widetilde{P}^{\lambda^{-1}(i)}
 \hookrightarrow G_{-1}\mathcal{C}_{\ell}(\lambda^{-1}(i))\cap
 \mathcal{C}_{\ell}^{\varepsilon}(\lambda^{-1}(i)) \hookrightarrow
 G_{-1}\mathcal{C}_{\ell}(\lambda^{-1}(i))
  \hookrightarrow \mathcal{D}_{\ell}^1(\lambda^{-1}(i))
  \label{CompositionOfInclusions}
\end{equation}
together with the decomposition (\ref{DecompositionOfD}) and the
homotopy equivalence in Lemma \ref{FirstDeformation}.

\section{The Salvetti Complexes for the Braid Arrangements}
\label{Configuration}

Consider the root system $A_{k-1}$, which is a collection of vectors in
the Cartan subalgebra $\mathfrak{h}_k$ of $\mathfrak{sl}_{k}(\R)$. We
may regard
\[
 \mathfrak{h}_k = \{(t_1,\cdots,t_k)\in\R^{k} \mid t_1+\cdots+t_k = 0 \}
\]
and the action of the Weyl group, which is the symmetric group of $k$
letters $\Sigma_{k}$, on $\mathfrak{h}_k$ is generated by the
reflections with respect to the hyperplanes
\[
 L_{i,j} = \{\bm{x} = (x_1,\cdots,x_k)\in \R^{k} \mid x_i=x_j\}.
\]
The complement of the complexification of this arrangement
\[
 \C^{k} -\bigcup_{i\neq j} L_{i,j}\otimes {\C}
\]
is nothing but the configuration space of $k$ points in $\C$, $F(\C,k)$,
and, under the identification (\ref{PermutationRepresentation}), it is
$\Sigma_k$-equivariantly homotopy equivalent to the complement 
\[
 \mathfrak{h}_k\otimes {\C} - \bigcup_{i\neq j}
 (L_{i,j}\cap\mathfrak{h}_k)\otimes {\C}. 
\]
In fact, the projections
\begin{eqnarray*}
 p & : & \R^k \longrightarrow \mathfrak{h}_k \\
 p\otimes{\C} & : & \C^k \longrightarrow \mathfrak{h}_k\otimes {\C}
\end{eqnarray*}
are $\Sigma_k$-equivariant and the latter induces a
$\Sigma_k$-equivariant homotopy equivalence between $F(\C,k)$ and
$\mathfrak{h}_k\otimes{\C} - \cup_{i\neq j} (L_{i,j}\cap\mathfrak{h}_k)
\otimes {\C}$.

The arrangement $\{L_{i,j} \mid i\neq j\}$ in $\R^k$ is 
denoted by $\mathcal{B}_n$. Let us denote the induced arrangement in
$\mathfrak{h}_k$ by $\mathcal{A}_{k-1}$. The purpose of this section is to
study the Salvetti complex for this arrangement and compare the skeletal
filtration on it to the filtration on $\mathcal{C}_2(k)$.

Note that the chambers of the real complement
\[
 \R^k -\bigcup_{i\neq j} L_{i,j}
\]
are labeled by elements in $\Sigma_k$. For
$\sigma \in \Sigma_k$, define corresponding points by
\begin{eqnarray*}
 \tilde{w}(\sigma) & = & (\sigma(1),\cdots, \sigma(k)) \\
  w(\sigma) & =  & p(\sigma(1), \cdots, \sigma(k)).
\end{eqnarray*}
Then each of them belongs to the chamber labeled by
$\sigma$. We have
\begin{eqnarray*}
 \mathcal{V}^0(\mathcal{B}_{k}) & = & \{\tilde{w}(\sigma) \mid
  \sigma\in\Sigma_k\} \\
 \mathcal{V}^0(\mathcal{A}_{k-1}) & = & \{w(\sigma) \mid
  \sigma\in\Sigma_k\}.
\end{eqnarray*}

More generally the faces in the $s$-th stratification (codimension $s$
facets) in the arrangement $\mathcal{B}_k$ are labeled by a partition
of $\{1,\cdots, k\}$ into ordered $k-s$ nonempty subsets. Recall that we
regard a partition as a surjective map
\[
 \lambda : \{1,\cdots, k\} \longrightarrow \{1,\cdots,k-s\}
\]
for some $s$.

\begin{definition}
 Given a partition $\lambda \in \Pi_k$, define a face by
 \[
 F_{\lambda} = \{(x_1,\cdots,x_k)\in \R^k \mid x_i<x_j \text{ if }
 \lambda(i)< \lambda(j) \text{ and } x_i=x_j \text{ if }
 \lambda(i) = \lambda(j)\}.
 \]
\end{definition}

\begin{lemma}
 The faces of $\mathcal{B}_k$ are
 \[
  \mathcal{F}(\mathcal{B}_k) = \{F_{\lambda} \mid \lambda \in \Pi_k\}.
 \]
 For brevity, we denote
\[
  \mathfrak{h}_{\lambda} = \mathfrak{h}_k\cap F_{\lambda}.
\]
 Then the facets of $\mathcal{A}_{k-1}$ are 
 \[
 \mathcal{F}(\mathcal{A}_{k-1}) = \{\mathfrak{h}_{\lambda} \mid 
 \lambda \in \Pi_k\}.
 \]
\end{lemma}

We can also use this map $\lambda$ to define a point in each facet.

\begin{definition}
 For a partition $\lambda \in \Pi_k$ , define
\begin{eqnarray*}
 \tilde{w}(\lambda) & = & (\lambda(1), \cdots, \lambda(k)) \\
 w(\lambda) & = & p(\tilde{w}(\lambda)).
\end{eqnarray*}
 Then $\tilde{w}(\lambda)$ and $w(\lambda)$ belong to the facet
 $F_{\lambda}$ and
 $\mathfrak{h}_{\lambda}$, respectively. Thus we have chosen vertex
 sets as
 \begin{eqnarray*}
  \mathcal{V}(\mathcal{B}_k) & = & \{\tilde{w}(\lambda) \mid \lambda \in
   \Pi_k\} \\
  \mathcal{V}(\mathcal{A}_{k-1}) & = & \{w(\lambda) \mid \lambda \in
   \Pi_k\}.
 \end{eqnarray*}
\end{definition}

The ordering of $\Pi_n$ corresponding to the ordering of
$\mathcal{F}(\mathcal{B}_k)$ is the following.

\begin{lemma}
 Order $\Pi_k$ as follows:
 \[
  \mu \le \lambda \Longleftrightarrow \lambda \text{ is a subdivision of }
 \mu,
 \]
 in other words
 \begin{eqnarray*}
  \lambda(i) = \lambda(j) & \Longrightarrow & \mu(i) = \mu(j) \\
  \lambda(i) < \lambda(j) & \Longrightarrow & \mu(i) \le \mu(j).
 \end{eqnarray*}
 Then $\Pi_k$ is isomorphic as posets to $\mathcal{F}(\mathcal{B}_k)$,
 hence to $\mathcal{F}(\mathcal{A}_{k-1})$.
\end{lemma}

Before we go on to discuss general cases in detail, let us take a look
at the simplest case, $k=2$.

\begin{example}
 There are three partitions of $\{1,2\}$:
 \[
  \Pi_2 = \{(1|2), (2|1), (1,2)\}.
 \]
 The facets of the arrangement $\mathcal{B}_2$ are
 \begin{eqnarray*}
  F_{(1|2)} & = & \{(x_1,x_2)\in\R^2 \mid x_1<x_2\} \\
  F_{(2|1)} & = & \{(x_1,x_2)\in\R^2 \mid x_2<x_1\} \\
  F_{(1,2)} & = & \{(x_1,x_2)\in\R^2 \mid x_1=x_2\} = M_{1,2}.
 \end{eqnarray*}
 and the facets of $\mathcal{A}_1$ are
 \begin{eqnarray*}
  \mathfrak{h}_{(1|2)} & = & F_{(1|2)}\cap\mathfrak{h}_2 \\
  \mathfrak{h}_{(2|1)} & = & F_{(2|1)}\cap\mathfrak{h}_2 \\
  \mathfrak{h}_{(1,2)} & = &  F_{(1,2)}\cap\mathfrak{h}_2.
 \end{eqnarray*}
 For each facet, the assigned vertex is given by
 \begin{eqnarray*}
  w((1|2)) & = & p(1,2) =
   (-\smallfrac{1}{2}, \smallfrac{1}{2}) \\
  w((2|1)) & = & p(2,1) = 
  (\smallfrac{1}{2}, -\smallfrac{1}{2}) \\
  w((1,2)) & = & p(1,1) = (0,0),
 \end{eqnarray*}
respectively.

We have only $0$-chains and $1$-chains. $0$-chains are
\[
 \mathfrak{h}_{(1,2)}, \mathfrak{h}_{(1|2)}, \mathfrak{h}_{(2|1)}
\]
and $1$-chains are
\[
 \mathfrak{h}_{(1|2)}>\mathfrak{h}_{(1,2)},
 \mathfrak{h}_{(2|1)}>\mathfrak{h}_{(1,2)}.
\]
In order to find simplices, we need to determine the facet-chamber
pairing. By Lemma \ref{PairingProperties}, we have
\begin{eqnarray*}
 \mathfrak{h}_{(1|2)}\cdot \mathfrak{h}_{(1|2)} & = &
  \mathfrak{h}_{(1|2)} \\
 \mathfrak{h}_{(1|2)}\cdot \mathfrak{h}_{(2|1)} & = &
  \mathfrak{h}_{(1|2)} \\
 \mathfrak{h}_{(2|1)}\cdot \mathfrak{h}_{(1|2)} & = &
  \mathfrak{h}_{(2|1)} \\
 \mathfrak{h}_{(2|1)}\cdot \mathfrak{h}_{(2|1)} & = &
  \mathfrak{h}_{(2|1)}.
\end{eqnarray*}
Since $(\mathcal{A}_2)_{\supset |\mathfrak{h}_{(1,2)}|} = \mathcal{A}_2$,
the unique facet containing $\mathfrak{h}_{(1|2)}$ is
$\mathfrak{h}_{(1|2)}$ itself. Since $\mathfrak{h}_{(1,2)}\cdot
\mathfrak{h}_{(1|2)}$ is a chamber contained in this
$\mathfrak{h}_{(1|2)}$, we have
\[
 \mathfrak{h}_{(1,2)}\cdot \mathfrak{h}_{(1|2)} = \mathfrak{h}_{(1|2)}
\]
and by the same reason
\[
 \mathfrak{h}_{(1,2)}\cdot \mathfrak{h}_{(2|1)} = \mathfrak{h}_{(2|1)}.
\]

Thus the $0$-simplices, i.e.\ vertices of the Salvetti complex are
\begin{eqnarray*}
 v(\mathfrak{h}_{(1,2)}, \mathfrak{h}_{(1,2)}\cdot \mathfrak{h}_{(1|2)})
  & = & v(\mathfrak{h}_{(1,2)}, \mathfrak{h}_{(1|2)}) \\
 & = & w((1,2)) + i(w((1|2))-w((1,2))) \\  
 & = & iw((1|2)) \\
 v(\mathfrak{h}_{(1,2)}, \mathfrak{h}_{(1,2)}\cdot \mathfrak{h}_{(2|1)})
  & = & v(\mathfrak{h}_{(1,2)}, \mathfrak{h}_{(2|1)}) \\
 & = & w((1,2))+ i(w((2|1))-w((1,2))) \\ 
 & = & iw((2|1)) \\
 v(\mathfrak{h}_{(1|2)}, \mathfrak{h}_{(1|2)}\cdot
 \mathfrak{h}_{(1|2)}) & = & v(\mathfrak{h}_{(1|2)},
 \mathfrak{h}_{(1|2)}) \\
 & = & w((1|2))+ i(w((1|2))-w((1|2))) \\
 & = & w((1|2)) \\ 
 v(\mathfrak{h}_{(2|1)}, \mathfrak{h}_{(2|1)}\cdot
 \mathfrak{h}_{(2|1)}) & = & v(\mathfrak{h}_{(2|1)},
 \mathfrak{h}_{(2|1)}) \\
 & = & w((2|1))+ i(w((2|1))-w((2|1))) \\
 & = & w((2|1)).
\end{eqnarray*}
We also have
\begin{eqnarray*}
 v(\mathfrak{h}_{(1|2)}, \mathfrak{h}_{(1|2)}\cdot
 \mathfrak{h}_{(2|1)}) & = & v(\mathfrak{h}_{(1|2)},
 \mathfrak{h}_{(1|2)}) \\
 & = & w((1|2))+ i(w((1|2))-w((1|2))) \\
 & = & w((1|2)) \\ 
 v(\mathfrak{h}_{(2|1)}, \mathfrak{h}_{(2|1)}\cdot
 \mathfrak{h}_{(1|2)}) & = & v(\mathfrak{h}_{(2|1)},
 \mathfrak{h}_{(2|1)}) \\
 & = & w((2|1))+ i(w((2|1))-w((2|1))) \\
 & = & w((2|1)).
\end{eqnarray*}
$1$-simplices are given by
\begin{eqnarray*}
 s(\mathfrak{h}_{(1|2)},\mathfrak{h}_{(1,2)}; \mathfrak{h}_{(1|2)}) & = &
  \Conv(\{v(\mathfrak{h}_{(1|2)}, \mathfrak{h}_{(1|2)}\cdot
  \mathfrak{h}_{(1|2)}), v(\mathfrak{h}_{(1,2)}, \mathfrak{h}_{(1,2)}\cdot
  \mathfrak{h}_{(1|2)})\}) \\
 & = & \Conv(\{w((1|2)), iw((1|2))\}) \\
 s(\mathfrak{h}_{(1|2)},\mathfrak{h}_{(1,2)}; \mathfrak{h}_{(2|1)}) & = & 
  \Conv(\{v(\mathfrak{h}_{(1|2)}, \mathfrak{h}_{(1|2)}\cdot
  \mathfrak{h}_{(2|1)}), v(\mathfrak{h}_{(1,2)}, \mathfrak{h}_{(1,2)}\cdot
  \mathfrak{h}_{(2|1)})\}) \\
 & = & \Conv(\{w((1|2)),iw((2|1))\}) \\
 s(\mathfrak{h}_{(2|1)},\mathfrak{h}_{(1,2)}; \mathfrak{h}_{(1|2)}) & = & 
  \Conv(\{v(\mathfrak{h}_{(2|1)}, \mathfrak{h}_{(2|1)}\cdot
  \mathfrak{h}_{(1|2)}), v(\mathfrak{h}_{(1,2)}, \mathfrak{h}_{(1,2)}\cdot
  \mathfrak{h}_{(1|2)})\}) \\
 & = & \Conv(\{w((2|1)), iw((1|2))\}) \\
 s(\mathfrak{h}_{(2|1)},\mathfrak{h}_{(1,2)}; \mathfrak{h}_{(2|1)}) & = & 
  \Conv(\{v(\mathfrak{h}_{(2|1)}, \mathfrak{h}_{(2|1)}\cdot
  \mathfrak{h}_{(2|1)}), v(\mathfrak{h}_{(1,2)}, \mathfrak{h}_{(1,2)}\cdot
  \mathfrak{h}_{(2|1)})\}) \\
 & = & \Conv(\{w((2|1)), iw((2|1))\}).
\end{eqnarray*}
Thus the Salvetti complex $\Sal(\mathcal{A}_1)$ is a $1$-dimensional
simplicial complex with four vertices
$(\smallfrac{1}{2}, -\smallfrac{1}{2})$,
$(-\smallfrac{1}{2}, \smallfrac{1}{2})$,
$(\smallfrac{i}{2}, -\smallfrac{i}{2})$,
$(-\smallfrac{i}{2}, \smallfrac{i}{2})$ and and four edges
$[(-\smallfrac{1}{2}, \smallfrac{1}{2}),
 (-\smallfrac{i}{2}, \smallfrac{i}{2})]$, 
$[(-\smallfrac{i}{2}, \smallfrac{i}{2}),
 (\smallfrac{1}{2}, -\smallfrac{1}{2})]$, 
$[(\smallfrac{1}{2}, -\smallfrac{1}{2}),
 (\smallfrac{i}{2}, -\smallfrac{i}{2})]$, 
$[(\smallfrac{i}{2}, -\smallfrac{i}{2}),
 (-\smallfrac{i}{2}, \smallfrac{i}{2})]$. It is the boundary of a square
 in $\mathfrak{h}_2^{\C}$ as is illustrated in the following picture.
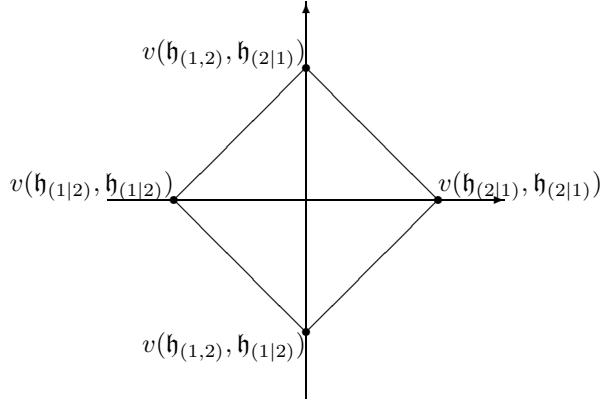
\begin{center}
 \begin{picture}(150,150)(-75,-75)
  \put(-75,0){\vector(1,0){150}}
  \put(0,-75){\vector(0,1){150}}
  \put(-50,0){\circle*{3}}
  \put(-50,0){\makebox(0,0)[br]{$v(\mathfrak{h}_{(1|2)}, \mathfrak{h}_{(1|2)})$}}
  \put(50,0){\circle*{3}}
  \put(50,0){\makebox(0,0)[bl]{$v(\mathfrak{h}_{(2|1)}, \mathfrak{h}_{(2|1)})$}}
  \put(0,-50){\circle*{3}}
  \put(0,-50){\makebox(0,0)[tr]{$v(\mathfrak{h}_{(1,2)}, \mathfrak{h}_{(1|2)})$}}
  \put(0,50){\circle*{3}}
  \put(0,50){\makebox(0,0)[br]{$v(\mathfrak{h}_{(1,2)}, \mathfrak{h}_{(2|1)})$}}
  \put(-50,0){\line(1,1){50}}
  \put(0,50){\line(1,-1){50}}
  \put(-50,0){\line(1,-1){50}}
  \put(0,-50){\line(1,1){50}}
 \end{picture}
\end{center}

 Note that under the composition of the inclusion maps in
 (\ref{CompositionOfInclusions}), the vertices $v(\mathfrak{h}_{(1|2)},
 \mathfrak{h}_{(1|2)})$, $v(\mathfrak{h}_{(1,2)}, \mathfrak{h}_{(2|1)})$,
 $v(\mathfrak{h}_{(2|1)}, \mathfrak{h}_{(2|1)})$, and
 $v(\mathfrak{h}_{(1,2)}, \mathfrak{h}_{(1|2)})$ correspond to cubes
 \begin{center}
  \begin{picture}(20,20)(0,0)
   \put(0, 0){\framebox(10,10){$1$}}
   \put(10,0){\framebox(10,10){$2$}}
  \end{picture},
  \begin{picture}(10,20)(0,0)
   \put(0, 0){\framebox(10,10){$2$}}
   \put(0,10){\framebox(10,10){$1$}}
  \end{picture},
  \begin{picture}(20,20)(0,0)
   \put(0, 0){\framebox(10,10){$1$}}
   \put(10,0){\framebox(10,10){$2$}}
  \end{picture},
  \begin{picture}(10,20)(0,0)
   \put(0, 0){\framebox(10,10){$1$}}
   \put(0,10){\framebox(10,10){$2$}}
  \end{picture},
 \end{center}
 respectively. And we have the following picture
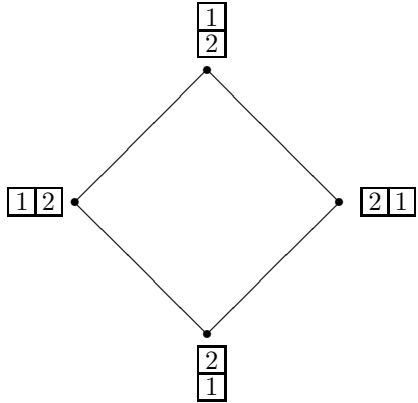
\begin{center}
 \begin{picture}(150,150)(-75,-75)
  \put(-50,0){\circle*{3}}
  \put(-55,0){\makebox(0,0)[r]{
  \begin{picture}(20,10)(0,0)
   \put(0, 0){\framebox(10,10){$1$}}
   \put(10,0){\framebox(10,10){$2$}}
  \end{picture}}}
  \put(50,0){\circle*{3}}
  \put(55,0){\makebox(0,0)[l]{
  \begin{picture}(20,10)(0,0)
   \put(0, 0){\framebox(10,10){$2$}}
   \put(10,0){\framebox(10,10){$1$}}
  \end{picture}}}
  \put(0,-50){\circle*{3}}
  \put(0,-55){\makebox(0,0)[t]{
  \begin{picture}(10,20)(0,0)
   \put(0, 0){\framebox(10,10){$1$}}
   \put(0,10){\framebox(10,10){$2$}}
  \end{picture}}}
  \put(0,50){\circle*{3}}
  \put(0,55){\makebox(0,0)[b]{
  \begin{picture}(10,20)(0,0)
   \put(0, 0){\framebox(10,10){$2$}}
   \put(0,10){\framebox(10,10){$1$}}
  \end{picture}}}
  \put(-50,0){\line(1,1){50}}
  \put(0,50){\line(1,-1){50}}
  \put(-50,0){\line(1,-1){50}}
  \put(0,-50){\line(1,1){50}}
 \end{picture}
\end{center}

 Let us consider the cell structure of $\Sal(\mathcal{A}_1)$. According
 to Lemma \ref{D(F,C)}, $1$-cells are
 \begin{eqnarray*}
  D(\mathfrak{h}_{(1,2)},h_{(1|2)}) & = &
   s(\mathfrak{h}_{(1|2)},\mathfrak{h}_{(1,2)};\mathfrak{h}_{(1|2)})
   \cup
   s(\mathfrak{h}_{(2|1)},\mathfrak{h}_{(1,2)};\mathfrak{h}_{(1|2)})\\
  & = & \Conv(\{w((1|2)), iw((1|2))\}) \cup \Conv(\{w((2|1)),
   iw((1|2))\}) \\
  & = &  \parbox{150pt}{\begin{picture}(150,75)(-75,-75)
  \put(-50,0){\circle*{3}}
  \put(50,0){\circle*{3}}
  \put(0,-50){\circle*{3}}
  \put(0,-55){\makebox(0,0)[t]{
  \begin{picture}(10,20)(0,0)
   \put(0, 0){\framebox(10,10){$1$}}
   \put(0,10){\framebox(10,10){$2$}}
  \end{picture}}}
  \put(-50,0){\line(1,-1){50}}
  \put(0,-50){\line(1,1){50}}
 \end{picture}} \\
  D(\mathfrak{h}_{(1,2)},h_{(2|1)}) & = &
   s(\mathfrak{h}_{(1|2)},\mathfrak{h}_{(1,2)};\mathfrak{h}_{(2|1)})
   \cup
   s(\mathfrak{h}_{(2|1)},\mathfrak{h}_{(1,2)};\mathfrak{h}_{(2|1)}) \\
  & = & \Conv(\{w((1|2)),iw((2|1))\}) \cup \Conv(\{w((2|1)),
   iw((2|1))\}) \\
  & = & \parbox{150pt}
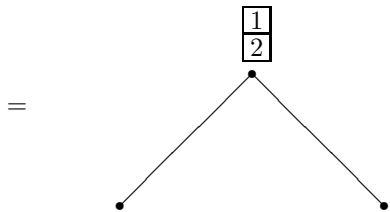
{ \begin{picture}(150,75)(-75,0)
  \put(-50,0){\circle*{3}}
  \put(50,0){\circle*{3}}
  \put(0,50){\circle*{3}}
  \put(0,55){\makebox(0,0)[b]{
  \begin{picture}(10,20)(0,0)
   \put(0, 0){\framebox(10,10){$2$}}
   \put(0,10){\framebox(10,10){$1$}}
  \end{picture}}}
  \put(-50,0){\line(1,1){50}}
  \put(0,50){\line(1,-1){50}}
 \end{picture}}
 \end{eqnarray*}
 and $0$-cells are
 \begin{eqnarray*}
  D(\mathfrak{h}_{(1|2)},\mathfrak{h}_{(1|2)}) & = &
   v(\mathfrak{h}_{(1|2)},\mathfrak{h}_{(1|2)}) \\
  & = & w((1|2)) \\
  & = & \parbox{20pt}{\begin{picture}(20,10)(0,0)
   \put(0, 0){\framebox(10,10){$1$}}
   \put(10,0){\framebox(10,10){$2$}}
  \end{picture}} \\
  D(\mathfrak{h}_{(2|1)},\mathfrak{h}_{(2|1)}) & = &
   v(\mathfrak{h}_{(2|1)},\mathfrak{h}_{(2|1)}) \\
  & = & w((2|1)) \\
  & = & \parbox{22pt}{\begin{picture}(20,10)(0,0)
   \put(0, 0){\framebox(10,10){$2$}}
   \put(10,0){\framebox(10,10){$1$}}
  \end{picture}}.
 \end{eqnarray*}

 Note that $1$-cells are labeled by cubes
 \begin{picture}(14,20)(-2,0)
   \put(0, 0){\framebox(10,10){$1$}}
   \put(0,10){\framebox(10,10){$2$}}
 \end{picture}
 and
 \begin{picture}(14,20)(-2,0)
   \put(0, 0){\framebox(10,10){$2$}}
   \put(0,10){\framebox(10,10){$1$}}
 \end{picture} in
 $F_{-1}\mathcal{C}_2(2)-F_{-2}\mathcal{C}_2(2)$ and the $0$-cells are
 labeled by cubes 
 \begin{picture}(24,10)(-2,0)
   \put(0, 0){\framebox(10,10){$1$}}
   \put(10,0){\framebox(10,10){$2$}}
  \end{picture}
 and
\begin{picture}(24,10)(-2,0)
   \put(0, 0){\framebox(10,10){$2$}}
   \put(10,0){\framebox(10,10){$1$}}
  \end{picture}
 in $F_{-2}\mathcal{C}_2(2)$. Thus the gravity
 filtration induces a filtration on the Salvetti 
 complex for $\mathcal{A}_1$ which coincides with the skeletal
 filtration up to the shift of filtration by $2$.
\end{example}

Let us return to the general case. For
$\lambda \in \Pi_k$ and $\sigma\in\Sigma_k$, the real parts of the
coordinates in the point $v(\lambda,\sigma)$ are determined by
$\lambda$ and the imaginary parts are determined by $\sigma$. To be more
precise, it is convenient to use the following symbols.

\begin{definition}
 \label{CubicalSymbol}
 For a partition $\lambda \in \Pi_k$ of rank $r$ and
 $\sigma \in \Sigma_k$, define a symbol $S(\lambda,\sigma)$ as follows:
 \begin{enumerate}
  \item For each $1\le i\le k-r$, draw vertically stacked squares $S_i$
	of length $|\lambda^{-1}(i)|$.
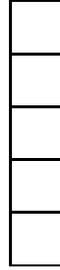
	\begin{center}
	 \begin{picture}(20,100)(0,0)
	  \put(0,0){\framebox(20,20)}
	  \put(0,20){\framebox(20,20)}
	  \put(0,40){\framebox(20,20)}
	  \put(0,60){\framebox(20,20)}
	  \put(0,80){\framebox(20,20)}
	 \end{picture}
	\end{center}
  \item Order $\lambda^{-1}(i)$ according to $\sigma$ and label each
	square in $S_i$  from bottom to top by elements in
	$\lambda^{-1}(i)$. For example, when
	$\lambda^{-1}(i) = \{i_1, i_2, i_3, i_4, i_5\}$ and if these
	numbers appear in $(\sigma(1),\cdots,\sigma(k))$ in the order
	\[
	 i_1,i_2,i_3,i_4,i_5
	\]
	then $S_i$ is labeled as
	
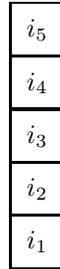
\begin{center}
	 \begin{picture}(20,100)(0,0)
	  \put(0, 0){\framebox(20,20){$i_1$}}
	  \put(0,20){\framebox(20,20){$i_2$}}
	  \put(0,40){\framebox(20,20){$i_3$}}
	  \put(0,60){\framebox(20,20){$i_4$}}
	  \put(0,80){\framebox(20,20){$i_5$}}
	 \end{picture}
	\end{center}
  \item Place $S_1, \cdots, S_{k-r}$ side by side from left to
	right. $S(\lambda,\sigma)$ is the resulting picture.
	
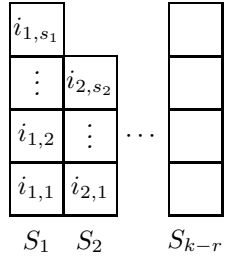
\begin{center}
	 \begin{picture}(80,100)(0,-20)
	  \put(0,-20){\makebox(20,20){$S_1$}}
	  \put(0, 0){\framebox(20,20){$i_{1,1}$}}
	  \put(0,20){\framebox(20,20){$i_{1,2}$}}
	  \put(0,40){\framebox(20,20){}}
	  \put(10,53){\makebox(0,0){$\vdots$}}
	  \put(0,60){\framebox(20,20){$i_{1,s_1}$}}

	  \put(20,-20){\makebox(20,20){$S_2$}}
	  \put(20, 0){\framebox(20,20){$i_{2,1}$}}
	  \put(20,20){\framebox(20,20){}}
	  \put(30,33){\makebox(0,0){$\vdots$}}
	  \put(20,40){\framebox(20,20){$i_{2,s_2}$}}

	  \put(40,20){\makebox(20,20){$\cdots$}}

	  \put(60,-20){\makebox(20,20){$S_{k-r}$}}
	  \put(60, 0){\framebox(20,20){}}
	  \put(60,20){\framebox(20,20){}}
	  \put(60,40){\framebox(20,20){}}
	  \put(60,60){\framebox(20,20){}}
	 \end{picture}
	\end{center}
 \end{enumerate}
\end{definition}

These are the symbols used in the picture of the hexagon in
Introduction.

The vertices of the simplicial complex $\Sal(\mathcal{A}_{k-1})$ are
given by 
\[
 v(\mathfrak{h}_{\lambda}, \mathfrak{h}_{\sigma}) = w(\lambda) +
 i(w(\sigma)-w(\lambda)) 
\]
for a partition $\lambda$ and a permutation $\sigma \in \Sigma_k$ which
is a subdivision of $\lambda$.

\begin{lemma}
 \label{VertexSymbols}
 There is a bijection between the set of vertices
 $\sk_0\Sal(\mathcal{A}_{k-1})$ and the set of symbols
 $\{S(\lambda,\sigma) \mid \lambda \in \Pi_k, \sigma\in\Sigma_k, 
 \lambda\le \sigma\}$.
\end{lemma}

\begin{definition}
 Define a filtration on the vertex set $\sk_0\Sal(\mathcal{A}_{k-1})$ by
 the number of distinct real coordinates:
 \[
  F_{-s}\sk_0\Sal(\mathcal{A}_{k-1}) = \{w(\lambda) +
 i(w(\sigma)-w(\lambda)) \mid |\Ima\lambda| \ge s\}.
 \]
 The associated filtration on the Salvetti complex by subcomplexes is
 denoted by
\begin{multline*}
  \emptyset = F_{-k-1}\Sal(\mathcal{A}_{k-1}) \subset
  F_{-k}\Sal(\mathcal{A}_{k-1}) \subset \cdots \subset \\
  F_{-s}\Sal(\mathcal{A}_{k-1}) \subset F_{-s+1}\Sal(\mathcal{A}_{k-1})
 \subset 
 \cdots \subset F_{-1}\Sal(\mathcal{A}_{k-1}) = \Sal(\mathcal{A}_{k-1}).
\end{multline*}
\end{definition}

Note that this is a filtration by $\Sigma_k$-subcomplexes.

\begin{lemma}
 The composition of the standard homotopy equivalences
 \[
  \Sal(\mathcal{A}_{k-1}) \hookrightarrow
 \mathfrak{h}_k\otimes {\C}-\bigcup_{i\neq j} ( L_{i,j}\cap
 \mathfrak{h}_k)\otimes {\C}
 \hookrightarrow \C^k-\bigcup_{i\neq j} L_{i,j}\otimes{\C} = F(\C,k)
 \longrightarrow \mathcal{C}_2(k)
 \]
 preserves the filtrations. Furthermore, the vertex $v(\lambda,\sigma)$
 is mapped to the symbol $S(\lambda,\sigma)$ which is regarded as an
 element of $\mathcal{C}_2(k)$ in an obvious way.
\end{lemma}

Let us denote this composition by
\[
 \varphi_k : \Sal(\mathcal{A}_{k-1}) \longrightarrow \mathcal{C}_2(k).
\]
Let us take a look at the induced map on the subquotients
\[
 F_{-s}\Sal(\mathcal{A}_{k-1})/F_{-s-1}\Sal(\mathcal{A}_{k-1})
 \longrightarrow  F_{-s}\mathcal{C}_2(k)/F_{-s-1}\mathcal{C}_2(k).
\]

Since the space of little $1$-cubes, $\mathcal{C}_1(k)$, is
$\Sigma_k$-equivariantly homotopy equivalent to $\Sigma_k$, Proposition
\ref{SubquotientsForC} give us the following description for the
subquotients for $\mathcal{C}_2(k)$.
\[
 F_{-s}\mathcal{C}_2(k)/F_{-s-1}\mathcal{C}_2(k)
 \simeq \bigvee_{\lambda \in \Pi_{k,k-s}}
 (\Sigma_{\lambda})_+ \wedge P^{\lambda}/\partial P^{\lambda} 
\]
where
\[
 \Sigma_{\lambda} = \Sigma_{\lambda^{-1}(1),\cdots,\lambda^{-1}(s)} =
 \Sigma_{\lambda^{-1}(1)}\times \cdots\times \Sigma_{\lambda^{-1}(s)}.
\]

The subquotients for the Salvetti complex can be easily found by
noticing that our filtration essentially coincides with the skeletal
filtration in Lemma \ref{D(F,C)}.

\begin{lemma}
 Let $\Sal(\mathcal{A}_{k-1})^{(s)}$ denote the $s$-skeleton of
 $\Sal(\mathcal{A}_{k-1})$ under the cell structure defined in Lemma
 \ref{D(F,C)}, then
 \[
  F_{-s}\Sal(\mathcal{A}_{k-1}) = \Sal(\mathcal{A}_{k-1})^{(k-s)}.
 \]
\end{lemma}

\begin{corollary}
 We have the following homeomorphism
 \begin{multline*}
  F_{-s}\Sal(\mathcal{A}_{k-1})/F_{-s-1}\Sal(\mathcal{A}_{k-1}) \cong \\
 \bigvee_{\lambda\in\Pi_{k,k-s}}
 \{\sigma \in \Sigma_k \mid \lambda\le \sigma\}_+ \wedge
 D(\lambda,\sigma)/\partial D(\lambda,\sigma).
 \end{multline*}
\end{corollary}

Note that, given a partition $\lambda$ of rank $k-s$, the set
$\{\sigma\in \Sigma_k \mid \lambda\le \sigma\}$ is in one-to-one
correspondence to the set
$\Sigma_{\lambda^{-1}(1)}\times \cdots\times \Sigma_{\lambda^{-1}(s)}$.
By investigation, we see that the map
$\varphi_k$ induces a $\Sigma_k$-equivariant homotopy equivalence on
each subquotient.

\begin{theorem}
 \label{Subquotient}
 $\varphi_k$ induces a $\Sigma_k$-equivariant homotopy equivalence for
 each $s$
 \[
 \varphi_k :
 F_{-s}\Sal(\mathcal{A}_{k-1})/F_{-s-1}\Sal(\mathcal{A}_{k-1})
 \rarrow{\Sigma_k} F_{-s}\mathcal{C}_2(k)/F_{-s-1}\mathcal{C}_2(k).
 \]
\end{theorem}

Recall that the spectral sequence (\ref{EMSS}) is induced from the
stable filtration on $\Omega^2\Sigma^2 X$ defined by the filtration
(\ref{gravity_filtration}) on little $2$-cubes. Thus we obtain the
following description of the $E^1$-term.

\begin{corollary}
 The $E^1$-term of the spectral sequence defined in \cite{Tamaki94} can
 be identified as follows
 \begin{eqnarray*}
  E^1_{-s,t} & = & \bigoplus_{k} E^1_{-s,t}(\mathcal{A}_{k-1}) \\
  & \cong &
 \bigoplus_{k}C_{k-s}(\Sal(\mathcal{A}_{k-1})) \otimes_{\Sigma_k}
 \tilde{h}_{t-k}(X^{\wedge k}) \\
  & \cong & \bigoplus_{k}C_{k-s}(\Sal(\mathcal{A}_{k-1}))
   \otimes_{\Sigma_k} \tilde{h}_{t}((\Sigma X)^{\wedge k}) \\
  & \cong & \bigoplus_k \Z\langle [D(\lambda,\sigma)] \mid \lambda \in
   \Pi_{k,k-s}, \sigma\in\Sigma_k, \lambda\le \sigma\rangle
   \otimes_{\Sigma_k} \tilde{h}_t((\Sigma X)^{\wedge k}). 
 \end{eqnarray*}
 This is an isomorphism of chain complexes.
\end{corollary}

Thus the $d^1$ of the spectral sequence is given by the cell structure
of $\Sal(\mathcal{A}_{k-1})$.
Lemma \ref{D(F,C)} also gives us a concrete description of the first
differential. Note that the cells in $\Sal(\mathcal{A}_{k-1})$ are
labeled by pairs of a partition and a permutation. In order to make an
explicit calculation, we use the following notations.

\begin{definition}
 \label{OrderedPartition}

  We denote a partition $\lambda$ of rank $r$ as a sequence of subsets
 \[
  \lambda = (\lambda^{-1}(1)|\cdots|\lambda^{-1}(k-r))
 \]
 or as a sequence of numbers separated by vertical lines
 \[
 \lambda = (i_{1,1},\cdots,i_{1,s_1}|\cdots|i_{k-j,1},\cdots,i_{k-r,s_{k-r}})
 \]
 when $\lambda^{-1}(1) = \{i_{1,1},\cdots,i_{1,s_1}\}$, $\cdots$,
 $\lambda^{-1}(k-r) = \{i_{k-r,1},\cdots,i_{k-r,s_{k-r}}\}$.

 The symmetric group $\Sigma_k$ acts on the set of partitions
 $\Pi_{k,r}$. As a representative of each orbit under this action, we
 may choose a partition which is an order-preserving map when regarded
 as a map
$\lambda : \{1,\cdots, k\} \longrightarrow \{1,\cdots,k-r\}$.
 Let $O_{k,r}$ be the set of order preserving
 partitions of rank $r$.
 
 An element $\lambda \in \Pi_{k,r}$ is said to be of type
 $(p_1,\cdots,p_{k-r})$, if 
 \[
  |\lambda^{-1}(1)| = p_1, \cdots, |\lambda^{-1}(k-r)| = p_{k-r}.
 \]
 The type of $\lambda$ is denoted by $t(\lambda)$.
\end{definition}

\begin{corollary}
 The $E^1$-term is given by
 \[
  E^{1}_{-s,t} \cong \bigoplus_n \Z\langle [D(\lambda,(1|\cdots|k))] \mid
   \lambda \in O_{k,k-s}\rangle \otimes \tilde{h}_t((\Sigma X))^{\otimes k}
 \]
\end{corollary}

\begin{example}
 Suppose $h_*(-)$ satisfies the K{\"u}nneth isomorphism and consider
 \begin{eqnarray*}
  d^1_{-1,t} : E^1_{-1,t} \longrightarrow E^1_{-2,t}.
 \end{eqnarray*}
 We have
 \begin{eqnarray*}
  E^1_{-1,*} & \cong & \bigoplus_{k}
   C_{k-1}(\Sal(\mathcal{A}_{k-1})) 
   \otimes_{\Sigma_k} \tilde{h}_{*}((\Sigma X))^{\otimes k} \\
  & \cong & \bigoplus_{k} \Z[\Sigma_k]\langle
   [D((1,\cdots,k),(1|\cdots|k))]\rangle 
   \otimes_{\Sigma_k} \tilde{h}_t((\Sigma X))^{\otimes k} \\
  & \cong & \bigoplus_{k} \Z\langle
   [D((1,\cdots,k),(1|\cdots|k))]\rangle \otimes \tilde{h}_t((\Sigma
   X))^{\otimes k} \\ 
  E^1_{-2,*} & \cong & \bigoplus_{k}
   C_{k-2}(\Sal(\mathcal{A}_{k-1})) \otimes_{\Sigma_k}
   \tilde{h}_*((\Sigma X))^{\otimes k} \\
  & \cong & \bigoplus_{k} \Z\langle [D(\lambda,(1|\cdots|k))] \mid
   \lambda \in O_{k,k-2}\rangle \otimes \tilde{h}_*((\Sigma X))^{\otimes k}.
 \end{eqnarray*}
 Consider the summand for $k=3$.
 Under the action of $\Sigma_3$, $\Pi_{3,1}$ has two orbits. One is
 represented by $(1|2,3)$ and the other is by $(1,2|3)$. The first
 differential is a map
 \begin{multline*}
  d^1_{-1,*} : \Z\langle([D(1,2,3),(1|2|3)])\rangle \otimes
  \tilde{h}_*(\Sigma X)^{\otimes 3} \\
  \longrightarrow \Z\langle
  [D(\lambda,(1|2|3))] \mid \lambda\in O_{3,1}\rangle \otimes
  \tilde{h}_*(\Sigma X)^{\otimes 3} \\
  =\Z\langle
  [D((1|2,3),(1|2|3))], [D((1,2|3),(1|2|3))] \rangle \otimes
  \tilde{h}_*(\Sigma X)^{\otimes 3}.
 \end{multline*}

 For elements $x_1, x_2, x_3 \in \tilde{h}_*(X)$, let us write the
 element $(\Sigma x_1)\otimes(\Sigma x_2)\otimes (\Sigma x_3) \in
 \tilde{h}_*(\Sigma X)^{\otimes 3}$ by $[x_1|x_2|x_3]$. Then we have
 \begin{multline*}
  d^1_{-1,t}([D((1,2,3),(1|2|3))]\otimes [x_1|x_2|x_3]) \\
  = \sum_{\lambda\in\Pi_{3,1}} [D(\lambda,
  \lambda\circ(1|2|3))]\otimes_{\Sigma_3} [x_1|x_2|x_3]
 \end{multline*}
 in $C_1(\Sal(\mathcal{A}_2))\otimes_{\Sigma_3}\otimes
 \tilde{h}_*(\Sigma X)^{\otimes 3}$, where each
 $D(\lambda,\lambda\circ(1|2|3))$ has the orientation induced from
 $D((1,2,3),(1|2|3))$.
 
 For
 $\lambda = (i_1|i_2,i_3) \in \Pi_{3,1}$ with $i_2<i_3$, let $\sigma =
 (i_1|i_2|i_3) \in \Pi_{3,0} = \Sigma_3$ then
 \[
  \lambda = (1|2,3)\circ\sigma.
 \]
 The difference of orientations between $D((1|2,3),(1|2|3))$ and
 $D((i_1|i_2,i_3),(1|2|3))$ is $\sgn(\sigma)$ where $\sgn$ is the sign
 function
 \[
 \sgn : \Sigma_3 \longrightarrow \{\pm 1\}.
 \]
 Then we have
 \begin{eqnarray*}
  [D(\lambda, \lambda\circ(1|2|3))]\otimes_{\Sigma_3}
   [x_1|x_2|x_3] 
  & = & \sgn(\sigma)[D((1|2,3), (1|2|3))]\cdot\sigma \otimes_{\Sigma_3}
   [x_1|x_2|x_3] \\
  & = & \sgn(\sigma)[D((1|2,3), (1|2|3))]\otimes_{\Sigma_3}
   [x_1|x_2|x_3]\cdot\sigma
 \end{eqnarray*}
 and
 \begin{multline*}
   d^1_{-1,t}([D((1,2,3),(1|2|3)]\otimes [x_1|x_2|x_3]) \\
   = [D((1|2,3), (1|2|3))]\otimes
   \left(\sum_{\sigma\in S_{1,2}} \sgn(\sigma)
    [x_1|x_2|x_3]\cdot\sigma\right) \\
   + [D((1,2|3), (1|2|3))]\otimes
   \left(\sum_{\sigma\in S_{2,1}} \sgn(\sigma)
    [x_1|x_2|x_3]\cdot\sigma\right)
 \end{multline*}
 where $S_{1,2}$ and $S_{2,1}$ are the set of $(1,2)$- and
 $(2,1)$-shuffles.

 For general $n$, $d^1_{-1,*}$ can be also described by $(p,q)$-shuffles.
\end{example}

In order to give a precise description of $d^1$, we need to compare
orientations of cells in $\Sal(\mathcal{A}_{k-1})$. Note that the action
of $\Sigma_k$ on $\Sal(\mathcal{A}_{k-1})$ is induced by permutation
of coordinates. For an $s$-cell $D(\lambda,\rho)$, choose an orientation.
Each $(s-1)$-cell in the boundary is assigned the orientation induced
from $D(\lambda,\rho)$. For two $(s-1)$-cells in the boundary,
$D(\mu,\mu\circ\rho)$ and $D(\mu',\mu'\circ\rho)$ if we write
\[
 D(\mu',\mu'\circ\rho) = D(\mu,\mu\circ\rho)\cdot\sigma
\]
with $\sigma\in \Sigma_k$, the $(s-1)$-cell $D(\mu',\mu'\circ\rho)$ has
two orientations, one induced from $D(\lambda,\rho)$ and the other
induced from $D(\mu,\mu\circ\rho)$. The difference of the orientations
is the determinant of the linear action of $\sigma$, $\sgn(\sigma)$.

We also need the following notation.

\begin{definition}
 An element $\sigma \in \Sigma_k$ is called a $(p_1,\cdots,p_s)$-shuffle,
 if $p_1+\cdots+p_s = k$ and 
 \begin{eqnarray*}
  & & \sigma(1)<\sigma(2)<\cdots<\sigma(p_1) \\
  & & \sigma(p_1+1) < \cdots < \sigma(p_1+p_2) \\
  & & \cdots \\
  & & \sigma(p_1+\cdots+p_{s-1}+1) < \cdots <
   \sigma(p_1+\cdots+p_s)=\sigma(k).
 \end{eqnarray*}
 The set of $(p_1,\cdots,p_s)$-shuffles is denoted by
 $S_{p_1,\cdots,p_s}$.
\end{definition}

\begin{corollary}
\label{d^1}
 For $x_1,\cdots, x_k \in \tilde{h}_*(X)$, let us denote
 \[
  [x_1|\cdots|x_k] = (\Sigma x_1)\otimes \cdots \otimes (\Sigma x_k) \in
 \tilde{h}_*(\Sigma X)^{\otimes k}.
 \]
 Then, for $\lambda \in O_{k,k-s}$, the
 first differential in the spectral sequence is given by 
 \begin{multline*}
  d^1_{-s,*}([D(\lambda,(1|\cdots|k))]\otimes[x_1|\cdots|x_k]) \\
  = \sum_{\tau \in O_{k,k-s-1}, \lambda<\tau} [D(\tau,(1|\cdots|k))]\otimes
 \left(\sum_{\rho\in S_{t(\tau)}}
 \sgn(\rho)[x_1|\cdots|x_k]\cdot\rho\right),
 \end{multline*}
 where $[D(\tau,(1|\cdots|k))]$ has the orientation induced from that of
 $[D(\lambda,(1|\cdots|k))]$.
 Recall that $S_{t(\tau)}$ is the set of shuffles of the same type as
 $\tau$.
\end{corollary}

It is worthwhile to note that the first differentials are
induced by space-level shuffles.
\[
 \bigvee_{\lambda\in O_{k,k-s}} S^{k-s}_{(\lambda,(1|\cdots|k))} \wedge
 (\Sigma X)^{\wedge k} 
 \longrightarrow \Sigma\left(\bigvee_{\tau\in O_{k,k-s-1, \lambda<\tau}}
 S^{n-s-1}_{(\tau,(1|\cdots|k))}
 \wedge (\Sigma X)^{\wedge k}\right),
\]
where $S^{k-s}_{(\lambda,(1|\cdots|k))}$ and
$S^{k-s-1}_{(\tau,(1|\cdots|k))}$ are copies of spheres $S^{k-s}$ and
$S^{k-s-1}$, respectively.

This map may be of some use to study maps between wedge powers of
suspended spaces.

\section{Higher Order Salvetti Complexes}
\label{higher_Salvetti}

Bj{\"o}rner and Ziegler briefly discussed
higher order oriented matroids in \S9.4 of \cite{Bjorner-Ziegler92}.
Let us construct higher order Salvetti complexes for oriented matroids
based on their idea. 

\subsection{Oriented Matroids and Real Arrangements}
\label{matroid_Salvetti}

In order to make the construction of the Salvetti complex described in
\S\ref{Complex} higher dimensional, we use oriented matroids. Let us
recall the relations between oriented matroids and real arrangements.

Let $\mathcal{A}=\{H_1,\cdots, H_n\}$ be a real central arrangement in a
vector space $V$. Choose a normal
vector $\bm{v}_i$ for each hyperplane $H_i$ and we obtain a vector
configuration $\mathcal{V} =\{\bm{v}_1,\cdots,\bm{v}_n\}$ in $V$.
Consider the set of linear dependencies of $\mathcal{V}$:
\begin{eqnarray*}
 \lindep(\mathcal{V}) & = & \{\bm{\lambda}=(\lambda_1,\cdots,\lambda_n)
  \in\R^n-\{\bm{0}\} \mid 
 \lambda_1\bm{v}_1+ \cdots + \lambda_n\bm{v}_n = \bm{0}\} \\
 & = & \{\bm{\lambda} : \mathcal{V} \longrightarrow \R \mid
  \bm{\lambda}(\bm{v}_1)\bm{v}_1+ \cdots +
  \bm{\lambda}(\bm{v}_n)\bm{v}_n = \bm{0}, \bm{\lambda}\neq 0\}.
\end{eqnarray*}

We say a linear dependency $\bm{\lambda}$ is minimal if it is no longer
a linear dependency when any one of $\lambda_i$'s is replaced with $0$.

By using the sign function
\[
 \sign : \R \longrightarrow \{+1,0,-1\}
\]
defined by
\[
 \sign(x) = \begin{cases}
	     +1, & x>0, \\
	     0, & x=0, \\
	     -1, & x<0,
	    \end{cases}
\]
we obtain a map
\[
 \sign\circ\bm{\lambda} : \mathcal{V} \longrightarrow \{+1,0,-1\}.
\]
The set
\[
 \mathcal{C}(\mathcal{V}) = \{\sign\circ\bm{\lambda} \mid
 \bm{\lambda} \in \lindep(\mathcal{V}):\text{ minimal}\}
\]
is called the set of signed circuits of $\mathcal{V}$. This collection
$\mathcal{C}(\mathcal{V})$ of signed circuits on $\mathcal{V}$ is
a typical example of an oriented matroid. In order to introduce a
general definition of oriented matroid, we regard $+1$, $0$, and $-1$
just as symbols and denote the set of three elements $\{+1,0,-1\}$ by
$S_1$. $S_1$ is equipped with a natural $\Z_2$-action. Note that a
function on a set $E$ with values in $S_1$
\[
 X : E \longrightarrow S_1
\]
can be regarded as a ``signed subset'' of $E$,
i.e.\ $X$ determines and is determined by two disjoint subsets of $E$
\begin{eqnarray*}
 X_{+} & = & X^{-1}(+1) \\
 X_{-} & = & X^{-1}(-1).
\end{eqnarray*}
The most fundamental signed subsets are the following Kronecker delta
functions $\delta_x^+$, $\delta_x^-$ defined for $x \in E$ by
\begin{eqnarray*}
 \delta_{x}^+(y) & = & \begin{cases}
			   +1, & x=y \\
			   0, & x\neq y
			  \end{cases} \\
 \delta_{x}^-(y) & = & \begin{cases}
			   -1, & x=y \\
			   0, & x\neq y.
			  \end{cases}
\end{eqnarray*}
A signed subset $X$ of $E$ can be regarded as a subset of
$E^{\pm} = \{\delta_{x}^+,\delta_{x}^{-} \mid x \in E\}$ satisfying
the following disjointness condition:
\[
 \delta_x^{\pm} \in X \Longrightarrow \delta_x^{\mp} \not\in X.
\]
For simplicity, we denote $\delta_x^{\pm}$ by $\pm x$.

\begin{remark}
 The adjoint
 \[
  \ad(\delta) : E \hookrightarrow \mathcal{P}(E)
 \]
 of the usual delta function
 \[
  \delta : E\times E \longrightarrow \{0,1\}
 \]
 embeds $E$ in $\Map(E,\{0,1\}) =\mathcal{P}(E)$
 as characteristic functions.

 Analogously, the signed delta functions
 \begin{eqnarray*}
  \delta^+ & : & E\times E \longrightarrow S_1 \\
  \delta^- & : & E\times E \longrightarrow S_1
 \end{eqnarray*}
 defined by
 \[
  \delta^{\pm}(x,y) = \begin{cases}
		       \pm 1, & x=y \\
		       0, & x\neq y
		      \end{cases}
 \]
 induce embeddings of $E$
 \begin{eqnarray*}
  \ad(\delta^+) & : & E \hookrightarrow \Map(E,S_1) \\
  \ad(\delta^-) & : & E \hookrightarrow \Map(E,S_1).
 \end{eqnarray*}
 The set $E^{\pm}$ is nothing but the union of the image of these maps
 \[
  E^{\pm} = \ad(\delta^+)(E) \cup \ad(\delta^-)(E).
 \]
\end{remark}

Note that the action of $\Z_2$ on the set $S_1$ induces an
action on $E^{\pm}$. 

\begin{definition}
 Given a set $E$, let
 \[
  \ast : E^{\pm} \longrightarrow E^{\pm}
 \]
 be the involution defined by $(\pm x)^{\ast} = \mp x$. A signed subset
 of $E$ is a subset $X$ of $E^{\pm}$ with
 \[
  X \cap X^{\ast} = \emptyset.
 \]
 The set of signed subsets of $E$ is denoted by
 $\mathcal{P}_{\pm}(E)$. 
\end{definition}

The following is a definition of oriented matroid in terms of subsets of
$E^{\pm}$ given in \cite{Gelfand-Rybnikov89}. This is essentially
identical to the definition by ``circuit axioms'' in
\cite{OrientedMatroids}.

\begin{definition}
 An oriented matroid on a set $E$ is a pair
 $\mathcal{M} = (E,\mathcal{C})$, where $\mathcal{C} \subset
 \mathcal{P}_{\pm}(E)-\{\emptyset\}$ is a collection of
 nonempty signed subset of $E$
 satisfying the following conditions: 
 \begin{enumerate}
  \item $X \in \mathcal{C}$ $\Longrightarrow$ $X^* \in \mathcal{C}$
  \item $X_1, X_2 \in \mathcal{C}$ and $X_1 \subset X_2\cup X_2^*$
	$\Longrightarrow$ $X_1=X_2$  or $X_1=X_2^*$
  \item $X_1, X_2 \in \mathcal{C}$, $e \in X_1\cap X_2^*$, and $X_1\neq
	X_2^*$ $\Longrightarrow$ there exists $Y \in \mathcal{C}$ such
	that
	\[
	 Y \subset (X_1\cup X_2)-\{e,e^*\}.
	\]
 \end{enumerate}
\end{definition}

\begin{example}
 Given a vector configuration
 $\mathcal{V} = \{\bm{v}_1,\cdots, \bm{v}_n\}$ in a real vector space
 $V$, the pair $\mathcal{M}(\mathcal{V}) = (\mathcal{V},
 \mathcal{C}(\mathcal{V}))$ is an oriented matroid.

 For a real hyperplane arrangement $\mathcal{A}$, the oriented matroid
 of the normal vector configuration $\mathcal{M}(\mathcal{V})$
 is independent of the choice of a normal vector configuration and is
 denoted by $\mathcal{M}(\mathcal{A})$. And the set of signed circuits is
 denoted by $\mathcal{C}(\mathcal{A})$.
\end{example}

Combinatorial structures of a real hyperplane arrangement
$\mathcal{A}$ can be described in terms of the associated oriented
matroid. For example, a face $F$ of $\mathcal{A}$ can be regarded as a
function
\[
 \tau_{F} : \mathcal{V}^{\pm} \longrightarrow S_1
\]
by
\[
 \tau_{F}(\bm{a}) = \begin{cases}
			  +1, & F \subset \{\bm{x} \mid \langle \bm{a},
		       \bm{x}\rangle >0\}, \\
			  0, & F \subset \{\bm{x} \mid \langle \bm{a},
		       \bm{x}\rangle =0\}, \\
			  -1, & F \subset \{\bm{x} \mid \langle \bm{a},
		       \bm{x}\rangle <0\}.
			 \end{cases}
\]
The following is a necessary and sufficient condition for such a
function to be associated with a face:

\begin{lemma}[Gel$'$fand-Rybnikov]
 Let $\mathcal{A}$ be a real central hyperplane arrangement and
 $\mathcal{V}$ be a vector configuration of normal vectors. For a
 function
 \[
  \tau : \mathcal{V}^{\pm} \longrightarrow
 S_1, 
 \]
 there exists a face $F$ with $\tau = \tau_F$ if and only if the
 following conditions hold:
 \begin{enumerate}
  \item $\tau(\pm\bm{a}) = \pm \tau(\bm{a})$, i.e.\ $\tau \in
	\Map_{\Z_2}(\mathcal{V}^{\pm},S_1)$, 
  \item For any signed circuit $X \in \mathcal{C}(\mathcal{A})$,
	\[
	 \tau(X) = \{0\} \text{ or } \tau(X) \supset \{+1,-1\}.
	\]
 \end{enumerate}
\end{lemma}

Thus we have a one-to-one correspondence between
$\mathcal{L}(\mathcal{A})$ and the set
 \[
   \{\tau \in
   \Map_{\Z_2}(\mathcal{V}^{\pm},S_1) \mid
   \tau(X) = \{0\} \text{ or } \tau(X) \supset \{+1,-1\} \text{ for }
   X\in\mathcal{C}(\mathcal{A})\}.
 \]
The poset structure on $\mathcal{L}(\mathcal{A})$ induces a poset
structure on this set.

\begin{lemma}
 For
 \[
  \varphi, \psi \in \{\tau \in
 \Map_{\Z_2}(\mathcal{V}^{\pm},S_1) \mid
 \tau(X) = \{0\} \text{ or } \tau(X) \supset \{+1,-1\} \text{ for }
 X\in\mathcal{C}\},
 \]
 $\varphi \le \psi$ if and only if
 \[
  \varphi(\bm{a}) \neq 0 \Longrightarrow \psi(\bm{a}) = \varphi(\bm{a}).
 \]
\end{lemma}

\begin{definition}
 For an oriented matroid $\mathcal{M} = (E,\mathcal{C})$,
 define
 \[
 \mathcal{L}(\mathcal{M}) = 
  \{\tau \in
 \Map_{\Z_2}(E^{\pm},S_1) \mid
 \tau(X) = \{0\} \text{ or } \tau(X) \supset \{+1,-1\} \text{ for }
 X\in\mathcal{C}\}.
 \]
 This is regarded as a poset by the following ordering:
 $\varphi \le \psi$ if and only if
 \[
  \varphi(\bm{a}) \neq 0 \Longrightarrow \psi(\bm{a}) = \varphi(\bm{a}).
 \]

 Elements in $\mathcal{L}(\mathcal{M})$ are called faces of
 $\mathcal{M}$ and a face $\varphi$ is called a chamber if
 $\varphi(\bm{a}) \neq 0$ for any $\bm{a} \in E^{\pm}$. The set of
 chambers is denoted by $\mathcal{L}^{(0)}(\mathcal{M})$.
\end{definition}

\begin{remark}
 For a real central arrangement $\mathcal{A}$ in a vector space $V$, we
 have a homotopy equivalence
 \[
  \mathcal{L}^{(0)}(\mathcal{A}) \simeq V-\bigcup_{L\in A} L
 \]
 if $\mathcal{L}^{(0)}(\mathcal{A})$ is regarded as a discrete
 $0$-dimensional complex.
\end{remark}

An important operation on faces of an oriented matroid is the following
matroid product, which corresponds to pairings in Proposition
\ref{matroid_product_of_faces}. 

\begin{definition}
 For $\varphi, \psi \in \Map(E,S_1)$, define
  $\varphi\circ\psi \in \Map(E,S_1)$
 by
 \begin{eqnarray*}
  (\varphi\circ\psi)(e) & = & \begin{cases}
			   \psi(e), & \varphi(e) \le \psi(e), \\
			   \varphi(e), & \textrm{otherwise}.
			  \end{cases} \\
  & = & \begin{cases}
	 \varphi(e), & \varphi(e) \neq 0, \\
	 \psi(e), & \varphi(e) = 0.
	\end{cases}
 \end{eqnarray*}
\end{definition}

\begin{lemma}
 The matroid product has the following properties: for $F_1, F_2\in
 \mathcal{L}(\mathcal{M})$, 
 \begin{enumerate}
  \item $F_1 \le F_1\circ F_2$,
  \item if $F_2\in \mathcal{L}^{(0)}(\mathcal{M})$ then $F_1\circ F_2
	\in \mathcal{L}^{(0)}(\mathcal{M})$,
  \item if $F_1\le F_2$ then $F_1\circ F_2 = F_2$,
  \item $F_1\circ 0 = F_1$.
 \end{enumerate}
\end{lemma}

It is well-known that the definition of oriented matroid can be given in
terms of faces (covectors). In order to describe the condition, it is
convenient to introduce the following.

\begin{definition}
 For $\sigma, \tau\in \Map_{\Z_2}(E^{\pm},S_1)$, define
 $S(\sigma, \tau) \subset E^{\pm}$ by
 \[
  S(\sigma,\tau) = \{x\in E^{\pm} \mid \sigma(x) = \tau(x)^* \neq 0\}.
 \]
 This is called the separation set of $\sigma$ and $\tau$.

 We say $\sigma$ is orthogonal to $\tau$ and write
 $\sigma\perp \tau$ if and only if $S(\sigma,\tau)$ and
 $S(\sigma, \tau^*)$ is both empty or both nonempty.
\end{definition}

The next proposition is called the ``covector axiom'' in
\cite{OrientedMatroids}.

\begin{proposition}
 A subset $\mathcal{L}$ of $\Map_{\Z_2}(E^{\pm},S_1)$ is the set
 of faces of an oriented matroid on $E$ if and only if
 it satisfies the following conditions:
 \begin{enumerate}
  \item $0 \in \mathcal{L}$,
  \item $\tau \in \mathcal{L}$ $\Longrightarrow$ $-\tau \in \mathcal{L}$,
  \item for $\tau_1, \tau_2 \in \mathcal{L}$, $\tau_1\circ \tau_2 \in
	\mathcal{L}$,
  \item for $\tau_1, \tau_2 \in \mathcal{L}$ and $x\in S(\tau_1, \tau_2)$,
	there exists $\tau_3 \in \mathcal{L}$ such that
	\begin{enumerate}
	 \item $\tau_3(x) = 0$
	 \item $\tau_3(y) = (\tau_1\circ \tau_2)(y) = (\tau_2\circ
	       \tau_1)(y)$ for all $y \not\in S(\tau_1,\tau_2)$.
	\end{enumerate}
 \end{enumerate}
\end{proposition}

In other words, a set of signed subsets $\mathcal{L}$ satisfying the
above conditions determines a set of signed circuits and vise
versa. Note that $\mathcal{P}_{\pm}(E)$ can be regarded as a subset of
$\Map_{\Z_2}(E^{\pm}, \{+1,0,-1\})$ by identifying $X \in
\mathcal{P}_{\pm}(E)$ with
\[
 \delta_X(x) = \begin{cases}
		1, & x\in X \\
		-1, & x\in -X \\
		0, & \textrm{otherwise}.
	       \end{cases}
\]
Under this identification, the definition of the set of faces is
given as follows:

\begin{proposition}
 Let $\mathcal{M} = (E,\mathcal{C})$ be an oriented matroid. Then the
 set of faces $\mathcal{L}$ is given by
 \[
  \mathcal{L} = \{\tau \in \Map_{\Z_2}(E^{\pm},S_1) \mid
 \tau\perp \varphi \textrm{ for all } \varphi\in\mathcal{C}\}.
 \]
\end{proposition}

\subsection{Higher Order Salvetti Complexes for Oriented Matroids}
\label{matroid_higher_Salvetti}

In this section we recall the definition of oriented $k$-matroid and the
$k$-dimension\-alization of an (ordinary) oriented matroid due to
Bj{\"o}rner and Ziegler \cite{Bjorner-Ziegler92}.

In order to understand their definitions, let us first recall the
definition of the Salvetti complex for an 
oriented matroid. The following description is due to Gel$'$fand and
Rybnikov \cite{Gelfand-Rybnikov89}.

\begin{definition}
 Let $S_2$ be the set of five elements $S_2 = \{0, e_1,-e_1,e_2,
 -e_2\}$. We define a partial order on $S_2$ by the following rule:
 \[
  0 <  \pm e_1 < \pm e_2.  
 \]

 For an oriented matroid $\mathcal{M} = (E,\mathcal{C})$, define
 $\mathcal{L}^{(1)}(\mathcal{M})$ by
 \[
 \mathcal{L}^{(1)}(\mathcal{M})  = 
   \{\varphi \in \Map_{\Z_2}(E^{\pm}, S_2-\{0\}) \mid
  \varphi(X) \subset \{\pm e_1\} 
  \text{ or } \varphi(X) \supset \{\pm e_2\} \text{ for }
  X\in\mathcal{C}\}.
 \]
\end{definition}

\begin{remark}
 The elements of the set $S_2$ are just symbols, but we may regard
 $\{e_1, e_2\}$ as the standard orthonormal basis of $\R^2$. $S_2$ is
 considered to be equipped with a $\Z_2$-action by changing signs.
\end{remark}

The following is an observation due to Gel$'$fand and Rybnikov, but we
include a proof in order to help the reader to understand an analogous
fact for the higher-order versions.

\begin{lemma}
 \label{L^(1)}
 The set $\mathcal{L}^{(1)}(\mathcal{M})$ is in one-to-one
 correspondence with the set of pairs
 \[
  \{(F,C) \in
 \mathcal{L}(\mathcal{M})\times
 \mathcal{L}^{(0)}(\mathcal{M}) \mid F\le C\}.
 \]
\end{lemma}

\begin{proof}
 In order to define a bijection, take $(F, C)$ in
 $\mathcal{L}(\mathcal{M})\times 
 \mathcal{L}^{(0)}(\mathcal{M})$ with $F\le C$. We use the following
 auxiliary $\Z_2$-equivariant maps
 \begin{eqnarray*}
  s_1 & : & \{\pm 1, 0\} \longrightarrow \{\pm e_1, 0\} \\
  s_2 & : & \{\pm 1, 0\} \longrightarrow \{\pm e_2, 0\}
 \end{eqnarray*}
 defined by
 \begin{eqnarray*}
  s_1(1) & = & e_1 \\
  s_1(0) & = & 0 \\
  s_2(1) & = & e_2 \\
  s_2(0) & = & 0
 \end{eqnarray*}
 For a face $F$, define $F\otimes e_i$ by
 \[
  F\otimes e_i = s_i\circ F : E \longrightarrow S_2.
 \]
 For maps
 \[
  F_1, F_2 : E \longrightarrow S_2
 \]
 we extend the matroid product as follows:
 \[
  (F_1\circ F_2)(x) = \begin{cases}
		       F_2(x), & F_2(x) > F_1(x) \\
		       F_1(x), & \textrm{otherwise.}
		      \end{cases}
 \]
  If $C$ is a chamber we obtain an element
 $(F\otimes e_1)\circ(C\otimes e_2) \in \mathcal{L}^{(1)}(\mathcal{M})$.

 In order to see this correspondence
 \[
  (F,C) \longmapsto (F\otimes e_1)\circ(C\otimes e_2)
 \]
 is a bijection, we need the following $\Z_2$-equivariant maps
 \begin{eqnarray*}
  \pi_1 & : & \{\pm e_1, \pm e_2\} \longrightarrow \{\pm 1, 0\} \\
  \pi_2 & : & \{\pm e_1, \pm e_2\} \longrightarrow \{\pm 1, 0\}
 \end{eqnarray*}
 defined by
 \begin{eqnarray*}
  \pi_1(e_1) & = & 1 \\
  \pi_1(e_2) & = & 1 \\
  \pi_2(e_1) & = & 1 \\
  \pi_2(e_2) & = & 0.
 \end{eqnarray*}
 Note that we have
 \begin{eqnarray*}
  \pi_1(F\otimes e_1) & = & F \\
  \pi_1(F\otimes e_2) & = & F \\
  \pi_2(F\otimes e_1) & = & F \\
  \pi_2(F\otimes e_2) & = & 0
 \end{eqnarray*}
 and 
 \begin{eqnarray*}
  \pi_1((F\otimes e_1)\circ (C\otimes e_2)) & = & \pi_1(F\otimes e_1)
   \circ \pi_1(C\otimes e_2) = F\circ C = C \\
  \pi_2((F\otimes e_1)\circ (C\otimes e_2)) & = & \pi_2(F\otimes e_1)\circ
   \pi_2(C\otimes e_2) = F\circ 0 = F.
 \end{eqnarray*}
 This completes the proof.
\end{proof}

\begin{definition}
 Define a partial order on $\mathcal{L}^{(1)}(\mathcal{M})$ as follows:
 $\varphi \le \psi$ if and only if
 \begin{eqnarray*}
  \varphi(\bm{a}) = e_0 & \Longrightarrow & \psi(\bm{a}) = e_0 \text{ or
   } \pm e_1 \\
  \varphi(\bm{a}) = -e_0 & \Longrightarrow & \psi(\bm{a}) = -e_0 \text{ or
   } \pm e_1 \\
  \varphi(\bm{a}) = e_1 & \Longrightarrow & \psi(\bm{a}) = e_1 \\
  \varphi(\bm{a}) = -e_1 & \Longrightarrow & \psi(\bm{a}) = -e_1.
 \end{eqnarray*}
 In other words, the ordering on $\mathcal{L}^{(1)}(\mathcal{M})$ is
 induced from the ordering on $S_2 = \{0, \pm e_0, \pm e_1\}$.
\end{definition}

\begin{definition}
 The first order Salvetti complex $\Sal^{(1)}(\mathcal{M})$ is the
 (geometric realization of the) order complex of
 $\mathcal{L}^{(1)}(\mathcal{M})$.
\end{definition}

As we have seen in the proof of Lemma \ref{L^(1)},
$\mathcal{L}^{(1)}(\mathcal{M})$ can be regarded as a subposet of
$(\mathcal{L}(\mathcal{M})\otimes e_1)\circ(\mathcal{L}(\mathcal{M})\otimes e_2)$.
The notion of oriented $2$-matroid was introduced by Bj{\"o}rner and
Ziegler in \cite{Bjorner-Ziegler92} by abstracting the properties of
this poset. At the end of their paper, they also introduced the notion
of oriented $\ell$-matroid for $\ell\ge 1$

We simplify their definition a little bit and introduce the notion of
symmetric oriented $k$-matroid.

\begin{definition}
 Let $S_{\ell}$ be the set
 \[
  S_{\ell} = \{0, e_1,-e_1,\cdots, e_{\ell},-e_{\ell}\}
 \]
 equipped with the following partial order:
 \begin{eqnarray*}
  0 & < & e_1, -e_1 \\
  e_1 & < & e_2, -e_2 \\
  -e_1 & < & e_2, -e_2 \\
  & \vdots & \\
  e_{\ell-1} & < & e_{\ell}, -e_{\ell} \\
  -e_{\ell-1} & < & e_{\ell}, -e_{\ell}
 \end{eqnarray*}
 and the obvious $\Z_2$-action.

 A signed $\ell$-vector on a set $E$ is a $\Z_2$-equivariant map
 \[
  \varphi : E^{\pm} \longrightarrow S_{\ell}.
 \]
 For signed $\ell$-vectors $\varphi$ and $\psi$, define the matroid product
 by
 \[
  (\varphi\circ\psi)(x) = \begin{cases}
			   \psi(x), & \psi(x) > \varphi(x) \\
			   \varphi(x), & \text{otherwise}
			  \end{cases}
 \]
 and the separation set by
 \[
  S(\varphi,\psi) = \{x\in E^{\pm} \mid \varphi(x)=\psi(x)^* \neq 0\}.
 \]
\end{definition}

Note that $S_{\ell}$ is equipped with a natural $\Z_2$-equivariant
$\Sigma_{\ell}$-action. This fact naturally leads us to the notion of
symmetric oriented $\ell$-matroid.

\begin{definition}
 A symmetric oriented $\ell$-matroid is a pair
 $\mathcal{M} = (E,\mathcal{L})$ of a set $E$ and
 $\mathcal{L} \subset \Map_{\Z_2}(E^{\pm}, S_{\ell})$ satisfying the
 following conditions: 
 \begin{enumerate}
  \item $0 \in \mathcal{L}$,
  \item $F \in \mathcal{L}$ $\Longrightarrow$ $-F \in \mathcal{L}$,
  \item $F \in \mathcal{L}$ $\Longrightarrow$
	$\sigma\cdot F \in \mathcal{L}$ for any $\sigma \in \Sigma_{\ell}$
  \item $F_1, F_2 \in \mathcal{L}$ $\Longrightarrow$
	$F_1\circ F_2 \in \mathcal{L}$,
  \item $F_1, F_2 \in \mathcal{L}$, $x \in S(F_1,F_2)$
	$\Longrightarrow$ $\exists F_3 \in \mathcal{L}$ such that
	$F_3(x) < F_1(x), F_2(x)$ and
	\[
	 F_3(y) = (F_1\circ F_2)(y) = (F_2\circ F_1)(y)
	\]
	for $y \not\in S(F_1,F_2)$.
 \end{enumerate}
\end{definition}

\begin{remark}
 A symmetric oriented $\ell$-matroid is a special case of an oriented
 $\ell$-matroid in the sense of Bj{\"o}rner and Ziegler. 
\end{remark}

There is a natural way to construct a symmetric oriented $\ell$-matroid
from an oriented ($1$-)matroid.

\begin{definition}
  For $1\le i\le \ell$, define a $\Z_2$-equivariant map
 \[
  s_i : \{\pm 1, 0\} \longrightarrow \{\pm e_i, 0\} \hookrightarrow
 S_{\ell} 
 \]
 by
 \begin{eqnarray*}
  s_i(1) & = & e_i \\
  s_i(0) & = & 0.
 \end{eqnarray*}
\end{definition}

\begin{lemma}
 Let $\mathcal{M}$ be an oriented ($1$-)matroid on a set $E$. 
 For a face $F \in \mathcal{L}(\mathcal{M})$ and $1\le i\le \ell$, define
 $F\otimes e_i \in \Map_{\Z_2}(E^{\pm}, S_{\ell})$ by
 \[
  (F\otimes e_i)(x) = s_i(F(x)).
 \]
 For simplicity, let us denote
 \[
  \mathcal{L}(\mathcal{M})\otimes\R^{\ell} = (\mathcal{L}(\mathcal{M})\otimes
 e_1)\circ \cdots \circ (\mathcal{L}(\mathcal{M})\otimes e_{\ell}).
 \]
 Then $(E, \mathcal{L}(\mathcal{M})\otimes\R^{\ell})$ is a symmetric
 oriented $\ell$-matroid.
\end{lemma}

\begin{proof}
 The first and the second conditions are obviously satisfied. The third
 and the fourth conditions follow from the fact that
 \[
  (F\otimes e_i)\circ (F'\otimes e_j) = (F'\otimes e_j)\circ (F\otimes
 e_i) 
 \]
 for $i\neq j$, and 
\[
  (F\otimes e_i)\circ (F'\otimes e_i) = (F\circ F')\otimes e_i.
\]
 The final condition follows from the corresponding condition for the
 covector axiom for oriented $1$-matroids.
\end{proof}

It is useful to identify elements of
$\mathcal{L}(\mathcal{M})\otimes\R^{\ell}$ with a sequence of faces of
$\mathcal{M}$.

\begin{lemma}
 For an oriented $1$-matroid $\mathcal{M}$, the set
 $\mathcal{L}(\mathcal{M})\otimes\R^{\ell}$ can be identified with  
 \[
  \{(F_1,\cdots, F_{\ell}) \in \mathcal{L}(M)^{\ell} \mid F_1\le \cdots
 \le F_{\ell}\}.
 \]
\end{lemma}

\begin{proof}
 For $0\le i\le \ell$, define a $\Z_2$-equivariant map
 \[
  \pi_i : S_{\ell} \longrightarrow \{\pm 1, 0\}
 \]
 by
 \begin{eqnarray*}
  \pi_i(e_j) & = & \begin{cases}
		    0, & j<i \\
		    1, & j\ge i.
		   \end{cases}
 \end{eqnarray*}
 Define a map
 \[
 c : \{(F_1,\cdots, F_{\ell}) \in \mathcal{L}(M)^{k} \mid F_1\le \cdots \le
 F_{\ell}\}  \longrightarrow
  \mathcal{L}(\mathcal{M})\otimes\R^{\ell}
 \]
 by
 \[
  c(F_1,\cdots,F_{\ell}) = (F_1\otimes e_1)\circ\cdots\circ
 (F_{\ell}\otimes e_{\ell}). 
 \]
 The inverse to this map is given by
 \[
  p(F) = (\pi_1\circ F,\cdots, \pi_{\ell}\circ F).
 \]
\end{proof}

\begin{definition}
 For a symmetric oriented $\ell$-matroid $\mathcal{M} = (E,\mathcal{L})$,
 define
 \[
  \mathcal{L}^{(0)} = \{\varphi \in \mathcal{L} \mid \varphi(x) \neq 0
 \text{ for all } x\in E\}
 \]
 An element of $\mathcal{L}^{(0)}$ is called a tope or a chamber.
\end{definition}

\begin{example}
 If $\mathcal{M}$ is an oriented $1$-matroid on a set $E$, then
 $(E, \mathcal{L}(\mathcal{M})\otimes\R^2)$
 is a symmetric oriented $2$-matroid. It is easy to see that
 \[
  \left(\mathcal{L}(\mathcal{M})\otimes\R^2\right)^{(0)} =
 \mathcal{L}^{(1)}(\mathcal{M}). 
 \]
\end{example}

\begin{definition}
 For an oriented $1$-matroid $\mathcal{M}$, let us denote
 \[
  \mathcal{L}^{(\ell)}(\mathcal{M}) =
 \left(\mathcal{L}(\mathcal{M})\otimes\R^{\ell+1}\right)^{(0)}.
 \]
 This set is regarded as a subposet of
 $\mathcal{L}(\mathcal{M})\otimes\R^{\ell+1}$, which has a poset
 structure induced from $S_{\ell+1}$.

 The $\ell$-th order Salvetti complex of $\mathcal{M}$ is the (geometric
 realization of the) order complex, or the classifying
 space of this poset
 \[
  \Sal^{(\ell)}(\mathcal{M}) =
 \left|\Delta(\mathcal{L}^{(\ell)}(\mathcal{M}))\right|.
 \]
\end{definition}

\begin{lemma}
 The set $\mathcal{L}^{(\ell)}(\mathcal{M})$ can be identified with
 \[
  \{(F_{\ell},\cdots, F_1,C) \in \mathcal{L}(M)^{\ell}\times
 \mathcal{L}^{(0)}(\mathcal{M}) \mid F_{\ell} \le \cdots \le F_1\le C\}.
 \]
\end{lemma}

Let $\mathcal{A}$ be a real central hyperplane arrangement in a vector
space $V$. The $\ell$-th order Salvetti complex for
$\mathcal{M}(\mathcal{A})$ is denoted by $\Sal^{(\ell)}(\mathcal{A})$.

$\Sal^{(\ell)}(\mathcal{A})$ can be embedded in
$V\otimes\R^{\ell+1}-\bigcup_{L\in\mathcal{A}} L\otimes\R^{\ell+1}$ as 
follows: we choose an interior point $v(F) \in V$ in each face $F$ of
the stratification of $V$ by $\mathcal{A}$ and define,
for $\varphi \in \mathcal{L}^{(k)}(\mathcal{A})$, 
\[
 v(\varphi) = v(\pi_0\circ\varphi)\otimes \bm{e}_0 +
 (v(\pi_0\circ\varphi)-v(\pi_1\circ\varphi))\otimes \bm{e}_1 + \cdots +
 (v(\pi_0\circ\varphi)-v(\pi_{\ell}\circ\varphi))\otimes \bm{e}_\ell,
\]
where $\{\bm{e}_0, \cdots, \bm{e}_{\ell}\}$ is the standard orthonormal
basis for $\R^{\ell+1}$.

$\Sal^{(\ell)}(\mathcal{M})$ is designed to extend the homotopy
equivalence  
\[
 \Sal^{(1)}(\mathcal{A}) \simeq
 V\otimes\R^2-\bigcup_{L\in\mathcal{A}}L\otimes\R^2
\]
for real central arrangements as follows.

\begin{theorem}
 \label{Embedding}
 The map
 \[
  v : \sk_0(\Sal^{(\ell)}(\mathcal{A})) = \mathcal{L}^{(\ell)}(\mathcal{A})
 \hookrightarrow V\otimes\R^{\ell+1} 
 \]
 induces an embedding
 \[
  v : \Sal^{(\ell)}(\mathcal{A}) \hookrightarrow
 \R^{\ell+1}-\bigcup_{L\in\mathcal{A}} L\otimes\R^{\ell+1}
 \]
 as a deformation retract.
\end{theorem}

This fact seems to be known to Bj{\"o}rner and Ziegler
\cite{Bjorner-Ziegler92}. The proof requires no novelty. The proof of
the analogous fact in \cite{Bjorner-Ziegler92} works if obvious
modifications are made. We include a sketch of proof in order to be
self-contained.

The idea of Bj{\"o}rner and Ziegler is to identify
$\Sal^{(k)}(\mathcal{A})$ as the (deformation retract of the) complement
of the link poset in $|\Delta(\mathcal{L}(\mathcal{A})\otimes
\R^{\ell+1}-\{0\})|$.

\begin{lemma}
 For a real essential central arrangement $\mathcal{A}$ in $V$, the cell
 complex $|\Delta(\mathcal{L}(\mathcal{A})\otimes \R^{\ell+1}-\{0\})|$ is
 homeomorphic to the unit sphere $S(V\otimes \R^{\ell+1})$.
\end{lemma}

\begin{definition}
 Define a subposet of $\mathcal{L}(\mathcal{A})\otimes \R^{\ell+1}-\{0\}$
 by
 \[
  \mathcal{K}_{\textrm{link}}(\mathcal{A}\otimes\R^{\ell+1}) = \{\varphi
 \in \mathcal{L}(\mathcal{A})\otimes \R^{\ell+1}-\{0\} \mid \varphi(e)=0
 \text{ for some } e\}.
 \]
 Namely this is the complement of $\mathcal{L}^{(\ell)}(\mathcal{A})$ in 
 $\mathcal{L}(\mathcal{A})\otimes \R^{\ell+1}-\{0\}$.
\end{definition}
 
The following is a well-known fact and can be found in
\cite{OrientedMatroids} as Lemma 4.7.27.

\begin{lemma}
 Let $P$ be a poset and $Q$ be a subposet. Then $|\Delta(Q)|$ is a
 deformation retract of $|\Delta(P)|-|\Delta(P-Q)|$.
\end{lemma}

The order complex of the poset
$\mathcal{K}_{\textrm{link}}(\mathcal{A}\otimes\R^{\ell+1})$ can be
identified with the intersection of $S(V\otimes\R^{\ell+1})$ and the
arrangement tensored with $\R^{\ell+1}$. 

\begin{lemma}
 We have the following homeomorphism:
 \[
  |\Delta(\mathcal{K}_{\textrm{link}}(\mathcal{A}\otimes\R^{\ell+1}))| \cong
 S(V\otimes\R^{\ell+1}) \cap
 \left(\bigcup_{L\in\mathcal{A}} L\otimes\R^{\ell+1}\right).
 \]
\end{lemma}

Thus $\Sal^{(\ell)}(\mathcal{A})$ is homotopy equivalent to
$V\otimes\R^{\ell+1}-\bigcup_{L\in\mathcal{A}} L\otimes\R^{\ell+1}$ and
Theorem \ref{Embedding} is proved.

$\Sal^{(\ell)}(\mathcal{A})$ is a simplicial complex as the order
complex of a poset. As we have done in \S\ref{Complex}, we can glue
simplices together to form a regular cell complex with fewer cells.
In stead of repeating the same argument, we use the following fact
observed in \cite{Bjorner-Ziegler92}. 

\begin{proposition}
 \label{cellular_realization}
 Let $K$ be a PL regular cell decomposition of a sphere $S^n$. For a
 subcomplex $L$, let $Q$ be the order ideal in the face poset $P$ of
 $K$. Then there exists a regular cell complex
 $L'$ with the following properties:
 \begin{enumerate}
  \item $L'$ is a subcomplex of the opposite regular cell decomposition
	of $K$,
  \item the face poset of $L'$ is isomorphic to $(P-Q)^{\op}$,
  \item $|L'|$ is homotopy equivalent to $|K|-|L|$.
 \end{enumerate}
\end{proposition}

Thus we can regard $\Sal^{(\ell)}(\mathcal{A})$ as a subcomplex
of $S(V\otimes\R^{\ell+1})$ under
a suitable regular cell decomposition of $S(V\otimes\R^{\ell+1})$. 

Now consider the braid arrangement $A_{k-1}$. We have
\[
 \Sal^{(\ell-1)}(\mathcal{A}_{k-1}) \simeq \mathfrak{h}_{k}\otimes
 \R^{\ell} - \bigcup_{1\le i<j\le k} L'_{i,j}\otimes \R^{\ell} \simeq
 F(\R^{\ell},k) 
 \simeq \mathcal{C}_{\ell}(k).
\]
The above homotopy equivalences are all $\Sigma_k$-equivariant.

The skeletal filtration on $\Sal^{(\ell-1)}(\mathcal{A}_{k-1})$ induces
the following spectral sequence.

\begin{theorem}
 \label{higher_gravity_spectral_sequence}
 There exists a spectral sequence for any homology theory
 \[
  E^1 \cong \bigoplus_k C_*(\Sal^{(\ell-1)}(\mathcal{A}_{k-1}))
 \otimes_{\Sigma_k}  
 \tilde{h}_*(X^{\wedge k}) \Longrightarrow
 h_*(\Omega^{\ell}\Sigma^{\ell} X),
 \]
 which is a direct sum of spectral sequences each of which strongly
 converges to the corresponding direct summand in
 (\ref{Snaith_splitting}). In particular, when $h_*(-)$ satisfies the
 strong form of K{\"u}nneth formula, the $E^1$-term is a functor of
 $h_*(X)$. 
\end{theorem}

Note that in the case $\ell=2$, the $E^1$-term is described as a functor
of $h_*(\Sigma X)$. The shift of degree by the suspension functor
$\Sigma$ is used to obtain an isomorphism between the $E^2$-term and
$\Cotor$. It is natural to expect that, when $\ell>2$, the $E^2$-term
can be also expressed as a certain derived functor by making
appropriate degree shifts in $h_*(X)$. 
\section{Concluding Remarks}
\label{remarks}

For the singular homology theory $H_*(-;k)$ with coefficients in a field
$k$, Smirnov \cite{Smirnov02-1} constructed a spectral sequence
converging to $H_*(\Omega^{\ell}Z;k)$ whose $E^1$-term is a functor of
$H_*(Z;k)$.

On the other hand, Ahearn and Kuhn studied the Goodwillie tower of the
functor $\Sigma^{\infty}\Map_*(K,Z)$, in particular
$\Sigma^{\infty}\Omega^{\ell} Z$, in
\cite{math/0109041,0806.3281} based on the analysis Arone did in
\cite{Arone99}. Kuhn remarks in \cite{0806.3281} that, when $K$ is a
sphere, the spectral sequence
\begin{equation}
 E^1(Z;\ell) \Longrightarrow h_*(\Omega^{\ell}Z) 
  \label{AAKSS}
\end{equation}
obtained from the Arone-Goodwillie tower must necessarily agree with
Smirnov's.

Arone's model describes the layers in the Arone-Goodwillie tower in
terms of little cubes. In particular, when $Z=\Sigma^{\ell} X$, the
tower coincides with the Snaith splitting and the spectral sequence
collapses at the $E^1$-term.

Our spectral sequence constructed in \S\ref{higher_Salvetti} is much
finer than the Arone-Ahearn-Kuhn spectral sequence (\ref{AAKSS}) for
$Z=\Omega^{\ell}X$. It is highly nontrivial. It would be interesting to
compare our spectral sequence and the Arone-Ahearn-Kuhn spectral
sequence in the following case 
\[
 E^1(\Omega^{\ell-m}\Sigma^{\ell} X;m) \Longrightarrow h_*(\Omega^{\ell}
 \Sigma^{\ell} X).
\]

\bibliographystyle{halpha}
\bibliography{%
\bibdir/mathAl,%
\bibdir/mathAn,%
\bibdir/mathA,%
\bibdir/mathB,%
\bibdir/mathBa,%
\bibdir/mathBe,%
\bibdir/mathBr,%
\bibdir/mathBu,%
\bibdir/mathCh,%
\bibdir/mathCo,%
\bibdir/mathC,%
\bibdir/mathDa,%
\bibdir/mathD,%
\bibdir/mathE,%
\bibdir/mathF,%
\bibdir/mathGa,%
\bibdir/mathGr,%
\bibdir/mathG,%
\bibdir/mathHa,%
\bibdir/mathHe,%
\bibdir/mathH,%
\bibdir/mathI,%
\bibdir/mathJ,%
\bibdir/mathKa,%
\bibdir/mathKo,%
\bibdir/mathK,%
\bibdir/mathLe,%
\bibdir/mathLu,%
\bibdir/mathL,%
\bibdir/mathMa,%
\bibdir/mathM,%
\bibdir/mathN,%
\bibdir/mathO,%
\bibdir/mathP,%
\bibdir/mathQ,%
\bibdir/mathR,%
\bibdir/mathSa,%
\bibdir/mathSe,%
\bibdir/mathS,%
\bibdir/mathT,%
\bibdir/mathU,%
\bibdir/mathV,%
\bibdir/mathW,%
\bibdir/mathX,%
\bibdir/mathY,%
\bibdir/mathZ,%
\bibdir/physics,%
\bibdir/personal}

\end{document}